\documentclass[]{article}

\usepackage{amssymb,latexsym,amsmath,amsthm}
\usepackage{graphicx,caption,subcaption,color}
\usepackage{stmaryrd}
\usepackage{multirow,multicol}
\usepackage{tikz}
\usetikzlibrary{shapes,arrows,decorations.pathreplacing}
\usetikzlibrary{positioning,fit}
\usepackage{import}
\DeclareMathOperator*{\supp}{supp}
\DeclareMathOperator*{\interior}{int}
\newtheorem{theorem}{Theorem}
\newtheorem{lemma}[theorem]{Lemma}

\newcommand{\norm}[1]{\left\|#1\right\|}

\newcommand{\R}{\mathbb{R}}

\newcommand{\dif}{\mathrm{d}}

\newcommand{\meansol}{\bar{p}}
\newcommand{\mssol}{p_{\text{ms}}}
\newcommand{\mesh}{\mathcal{T}}
\newcommand{\edge}{\mathcal{E}}
\newcommand{\msspace}{V_{\text{ms}}}
\newcommand{\aux}{V_{\text{aux}}}
\newcommand{\glo}{V_{\text{glo}}}
\newcommand{\osdomain}[1]{K_{#1,m}}
\newcommand{\auxbasis}[2]{\phi_{#1}^{(#2)}}
\newcommand{\msbasis}[2]{\psi_{#1,\text{ms}}^{(#2)}}
\newcommand{\globasis}[2]{\psi_{#1}^{(#2)}}
\newcommand{\kij}[2]{K_{#1}^{(#2)}}
\newcommand{\lij}{L_i^{(j)}}
\newcommand{\Rglo}{R_{\text{glo}}}
\newcommand{\Rms}{R_{\text{ms}}}

\newcommand{\lsquare}[1]{[L^2(#1)]^2}
\newcommand{\lweight}[2]{\mathbb{L}^2(#1;#2)}
\newcommand{\unit}[1]{\mathbf{e}_{#1}}

\graphicspath{{figures/}}

\begin{document}

\title{Analysis of Non-local Multicontinuum Upscaling for Dual Continuum Model}
\author{
Jingyan Zhang\thanks{Department of Mathematics, Texas A\&M University, College Station, TX 77843, USA (\texttt{jingyanzhang@math.tamu.edu})}
\and
Siu Wun Cheung\thanks{Department of Mathematics, Texas A\&M University, College Station, TX 77843, USA (\texttt{tonycsw2905@math.tamu.edu})}
}\maketitle

\begin{abstract}
In this paper, we develop and analyze a rigorous multiscale upscaling method for dual continuum model, 
which serves as a powerful tool in subsurface formation applications. 
Our proposed method is capable of identifying different continua and capturing non-local transfer and effective properties 
in the computational domain via constructing localized multiscale basis functions. 
The construction of the basis functions consists of solving local problems defined on oversampling computational region, 
subject to the energy minimizing constraints that the mean values of the local solution are zero in all continua except for the one targeted. 
The basis functions constructed are shown to have good approximation properties. 
It is shown that the method has a coarse mesh dependent convergence. 
We present some numerical examples to illustrate the performance of the proposed method. 

\end{abstract}

\section{Introduction}


Subsurface formations exist in a variety of practical applications. 
In reservoir simulation, the material properties within fractures and background media can be significantly distinct. 
In the case of large fractures, Discrete Fracture Model (DFM) and Embedded Frature Model (EFM) are used to explicitly define fracture networks with accuracy 
\cite{karimifard2001dfm,karimifard2004edfm,garipov2016dfm}. 
Explicitness for complex models is naturally in need of large system of equations, which will lead to large computational costs. 
Moreover, due to the multiple scales and high contrast properties intrinsic to the reservoir, performing high fidelity simulations for fractured porous media is a complex enough to cause extremely high computational costs. 

To overcome this difficulty, research effort had been devoted to construct numerical solvers 
on a coarse grid, which is typically much coarser than the fine grid which captures all the heterogeneities in the medium properties. 
Typical approaches involve computing upscaled effective properties in each local coarse-grid block or representative volume \cite{dur91,weh02}. 
Such approaches are known to be insufficient when more than one important modes exist in the same coarse block or representative volume.  
In these cases, more efficient upscaling methods, typically the multicontinuum models, 
are employed \cite{arbogast1990derivation,barenblatt1960basic,kazemi1976numerical,pruess1982fluid,warren1963behavior,wu1988multiple}. 
In such approaches, several effective properties are formulated in each coarse block and 
interaction terms are defined to characterize the transfer between different continua. 
Another class of methods are the multiscale methods, including
Heterogeneous Multiscale Methods (HMM) \cite{ee03,abdulle05,emz05}, 
Variational Multiscale Methods (VMS) \cite{hfmq98,hughes2007variational,Iliev_MMS_11,calo2011note} 
and Multiscale Finite Element Method (MsFEM) \cite{hw97,ehw99}.
Similar to upscaling approaches, multiscale scale is to construct numerical solvers 
on the coarse grid, which is typically much coarser than the fine grid
which captures all the heterogeneities in the medium properties. 
Instead of computing the effective medium properties, multiscale basis functions which are 
responsible for capturing the local oscillatory effects of the solution are constructed and coarse-scale macroscopic equations are formulated.
The solution of the coarse-scale system can then be used to recover 
fine-scale information with the mutliscale basis functions. 

However, for more complex high-contrast heterogeneous media, each local coarse region contains several high-conductivity regions and 
multiple multiscale basis functions are required to represent the local solution space. 
The aforementioned multiscale methods typically use one basis function per coarse region,
which is insufficient and may give rise to large error. To this end, it is crucial to systemically enrich the 
multiscale space with suitable fine-scale information for low-dimensional solution representation. 
One such approach is the Generalized Multiscale Finite Element Method (GMsFEM) 
\cite{egw10,egh12,chung2014adaptive,chung2016adaptive}, 
which involves the construction appropriate snapshots space which consists of fine-scale data for solution representation by local snapshot problems 
and the construction of multiscale basis functions by performing feature extraction through local spectral decompositions to the snapshot space. 
Since multiscale basis functions are identified for multiscale solution representation in high-contrast heterogenerous media,
multiple basis functions from the spectral problem are required to attain a small error.
By introducing adaptivity \cite{efendiev2016online,chung2016adaptive}, one can add multiscale basis functions in selected regions. 
The connection of GMsFEM to multicontinuum models is discussed in \cite{chung2017coupling}, 
and GMsFEM are successfully applied to multicontinuum models that originate from fracture models and 
contain nonlinearities \cite{wang2020vug,cheung2018mc,park2019mc,park2020richard}. 

More recently, a combination of GMsFEM and localization has been discussed in \cite{chung2018constraint} 
as the approach of Constraint Energy Minimizing GMsFEM (CEM-GMsFEM). 
The method uses oversampling computational regions for the construction of multiscale basis functions. 
The first step is to find the auxiliary multiscale basis functions by GMsFEM. 
The second step is to construct multiscale basis functions by minimizing energy functionals subject to certain constraints, the purpose of which is to localize the multiscale basis functions. 
The method has been applied to various discretization and model problems \cite{chung2018mixed,cheung2019dg,cheung2018mc,cheung2020wave}, and 
it has been theoretically and numerically verified that the multiscale solutions spanned by the multiscale basis functions
in CEM-GMsFEM have both spectral convergence and mesh dependent convergence. 

The construction of auxiliary multiscale functional space is key to identify high contrast channels and fractures in multicontinuum models.  
However, the adoption of GMsFEM in obtaining the auxiliary space is of relatively high computational cost. 
To modify the method, under the assumption that one knows the fracture network in each coarse-grid block, 
the method of nonlocal multicontinuum (NLMC) is proposed in \cite{chung2018nlmc,zhao2020nlmc}
and applied to fracture models \cite{maria2019nlmc,maria2019edfm}. 
The aforementioned assumption is a main drawback to NLMC but is a scenario common in practice. 
Instead of solving local spectral problems GMsFEM, the auxiliary basis functions are defined to represent 
coarse-scale solution average in a straightforward manner. 
As is in CEM-GMsFEM, the multiscale basis functions are solved from problems formulated to minimize the local energy in an oversampling domain. 
The mass transfers between fractures and matrix is therefore non-local. 

In this work, we develop and analyze the NLMC method for a dual continuum model. 
The auxiliary basis functions are simply defined in each coarse block for each continuum, 
which represents fractures and matrix. 
To be more precise, in each oversampling domain, the auxiliary basis functions are constant in a continuum, and 
each has mean value one for the chosen continuum and zero otherwise. 
Out of the oversampling domain, the value of the basis functions are zero. 
The degrees of freedom is the same as the number of the continua, which is the minimal number needed to represent the heterogeneous property of the reservoir. 
To obtain the multiscale basis functions, we solve local minimization problems in oversampling computational domain. 
We show that the minimizer has a good decay property. 
With a proper number of oversampling layers, the basis functions derived can well capture the fine-grid information. 
Moreover, we show that the method has a convergence dependent on coarse mesh size. 
We also present some numerical examples to depict the performance of the method.

The paper is organized as follows. 
In Section \ref{sec:model} we present the dual continuum model. 
The proposed method is introduced in Section \ref{sec:method} and analyzed in Section \ref{sec:analysis}. 
In Section \ref{sec:numerical} some numerical experiments are demonstrated to confirm the the theory. 
The paper ends with conclusions in Section \ref{sec:conclusions}.

\section{Dual continuum Model}\label{sec:model}

We consider the following dual continuum model \cite{barenblatt1960basic, douglas1990dual, warren1963behavior}
\begin{equation}
\begin{split}
c_1 \dfrac{\partial p_1}{\partial t} - \text{div}(\kappa_1 \nabla p_1) + \sigma(p_1 - p_2) = f_1, \\
c_2 \dfrac{\partial p_2}{\partial t} - \text{div}(\kappa_2 \nabla p_2) - \sigma(p_1 - p_2) = f_2,
\end{split}
\label{eq:dc}
\end{equation}
in a computational domain $\Omega \subset \mathbb{R}^{{2}}$.  
{Here, for $i = 1,2$, $c_i$ is the compressibility, $p_i$ is the pressure, 
$\kappa_i$ is the permeability, 
and $f_i$ is the source function for the $i$-th continuum. 
In addition, the continua are coupled through the mass exchange, 
and $\sigma$ is a parameter which accounts for the strength of mass transfer between the continua.  
One particular application of the dual continuum model \eqref{eq:dc} 
is to represent the global interactive effects of the unresolved fractures and the matrix. 
In this work, we consider high-contrast channelized media. 
We prescribe the initial condition $p_i(0,\cdot) = p_i^0$ in $\Omega$ and 
the boundary condition $p_i(t, \cdot) = 0$ on $\partial \Omega$ for $t > 0$. 

Let $V = [H^1_0(\Omega)]^2$. 
Also, for a subdomain $D \subset \Omega$, we denote the restriction of $V$ on $D$ by $V(D)$.
The weak formulation of \eqref{eq:dc} then reads:
find $p = (p_1, p_2)$ such that $p(t, \cdot) \in V$ and
\begin{equation}
\begin{split}
c \left(\dfrac{\partial p}{\partial t}, v \right) + a_Q(p, v) = (f,v),
\end{split}
\label{eq:sol_weak}
\end{equation}
for all $v = (v_1, v_2)$ with $v(t, \cdot) \in V$.
{Here, $(\cdot, \cdot)$ denotes the standard $L^2(\Omega)$ inner product. 
Moreover, }the bilinear forms are defined as:
\begin{equation}
\begin{split}
c_i ({p}_i, v_i) & = 
{ \int_{\Omega} c_i (x) {p}_i v_i \, \dif x,}\\
c({p}, v) & = \sum_i c_i ({p}_i,v_i), \\
a_i (p_i, v_i) & = 
{ \int_{\Omega} \kappa_i (x) \nabla p_i \cdot \nabla v_i \, \dif x},\\
a(p,v) & = \sum_i a_i(p_i, v_i), \\
q(p,v) & = \sum_i \sum_{i'}  \int_\Omega {\sigma} (p_i - p_{i'}) v_i \, \dif x,\\
a_Q(p,v) & = a(p,v) + q(p,v). 
\end{split}
\label{eq:bilinear}
\end{equation}

\section{Method description}\label{sec:method}
In this section, we will describe our proposed method in detail. 
To start with, we introduce the concepts of coarse and fine meshes. 
We start with a plain partition of calculation domain $\Omega$, $\mesh^H$. 
This partition is called a coarse mesh, which does not necessarily resolve any multiscale features. 
We denote one element in $\mesh^H$ as $K$ and name it as a coarse element. 
Here, $H>0$ is the coarse mesh size. 
We denote the number of coarse elements and coarse grid nodes as $N$ and $N_c$ respectively. 
The collection of all coarse element edges is called $\edge^H$. 
To sufficiently resolve the solution, we refine the coarse mesh $\mesh^H$ into a fine mesh $\mesh^h$, where $h>0$ is called the fine mesh size. 
We remark that the fine grid system is only used in locally solving process, where all local problems are solved continuously. 
Therefore, we don't consider fine grid in our analysis hereinafter. 
A demonstration of coarse and fine meshes is given in Figure~\ref{fig:mesh}. 

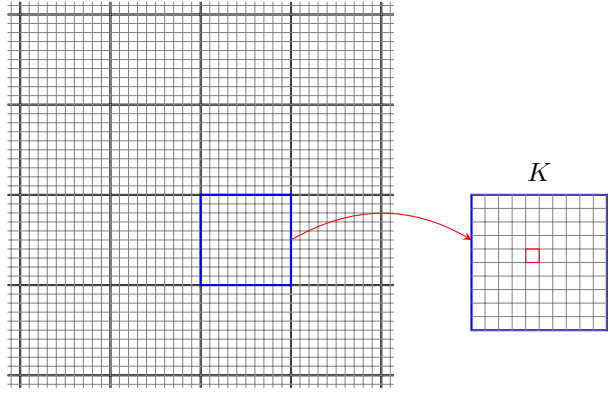
\begin{figure}[ht!]
\centering
\begin{tikzpicture}[scale = 1.2]
	\draw[step=1cm,black,thick] (-2.14,-2.14) grid (2.14,2.14);
	\draw[step=0.1cm, gray, very thin] (-2.14,-2.14) grid (2.14,2.14);
	\draw[blue,thick] (0,0) rectangle (1,-1);
	\path[red,-stealth] (1,-0.5) edge[bend left] (3,-0.5);
	\draw[blue,thick] (3,0) rectangle (4.5,-1.5);
	\draw[step=0.15cm, gray, very thin] (3,0) grid (4.5,-1.5);
	\draw[red,thin] (3.6,-0.6) rectangle (3.75, -0.75);
	\draw (3.75,0.25) node(K) [draw=white] {$K$};
	
\end{tikzpicture}

\caption{An illustration of the coarse and fine meshes.}
\label{fig:mesh}
\end{figure}

Next, we clarify more notions concerning every coarse element. 
Letting $K_j \in \mesh^H$ be the $j$-th coarse element, an oversampling domain $\osdomain{j}$ is defined by expanding $K_j$ with $m$ layers of coarse elements in $\Omega$. 
An illustration is given in Figure~\ref{fig:osdomain}. 
Moreover, similar to the partition of computational domain $\Omega$, for $i = 1,2$, we denote $K_j = \cup_{l=1}^{L_i^{(j)}} \kij{l}{i,j}$, where $L_i^{(j)}$ denotes the number of high conductive channels plus matrix within $K_j$ for continua $i$. 

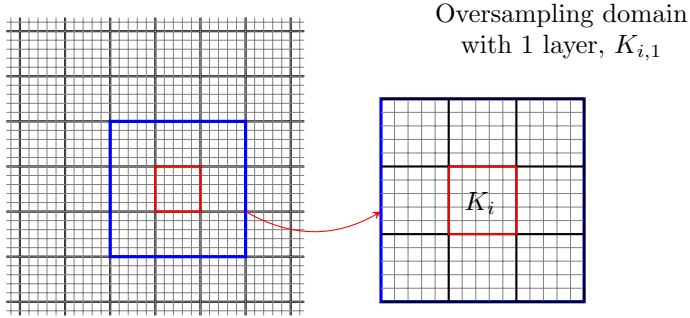
\begin{figure}[ht!]
\centering
\begin{tikzpicture}[scale = 1.2]
	\draw[step=0.5cm,black,thick] (-1.15,-1.65) grid (2.15,1.65);
	\draw[step=0.1cm, gray, very thin] (-1.15,-1.65) grid (2.15,1.65);
	\draw[blue,very thick] (0,-1) rectangle (1.5,0.5);
	\draw[red,thick] (0.5,-0.5) rectangle (1,0);
	
	\path[red,-stealth] (1.5,-0.5) edge[bend right] (3,-0.5);
	\draw[blue,very thick] (3,0.75) rectangle (5.25,-1.5);
	\draw[step=0.15cm, gray, very thin] (3,0.75) grid (5.25,-1.5);
	\draw[step=0.75cm, black,thick] (3,0.75) grid (5.25,-1.5);
	\draw[red,thick] (3.75,-0.75) rectangle (4.5, 0);
	\draw (5,1.5) node(x) [draw=white] {\begin{tabular}{c}Oversampling domain\\ with 1 layer, $K_{i,1}$\end{tabular}};
	
	\draw (4.1,-0.4) node(ki){$K_i$};
	
\end{tikzpicture}
\caption{An illustration of oversampling domain.}
\label{fig:osdomain}
\end{figure}

We now proceed to describing the proposed method step by step. 

\textbf{Step 1. Definition of auxiliary basis functions.} 
Within $K_j$($j = 1,\ldots, N$), for each $l = 1,\ldots, L_i^{(j)}$, we directly define our auxiliary basis function $\auxbasis{l}{i,j} \in P^0(\mesh^h(\osdomain{j}))$, where $\mesh^h(\osdomain{j})$ denotes the restriction of $\mesh^h$ on $\osdomain{j}$, as 
\begin{equation}
	\auxbasis{l}{i,j} = \dfrac{1}{|K_l^{(i,j)}|}\mathcal{I}_{\kij{l}{i,j}}.
	\label{eq:auxbasis}
\end{equation}
Here, 
$|\kij{l}{i,j}|$ denotes the area of $\kij{l}{i,j}$ and 
$\mathcal{I}_{\kij{l}{i,j}}$ is the characteristic function.  
The local auxiliary space for continua $i$ is constructed as
\begin{equation}
	\aux^{(i,j)} = \mathrm{span}\{\auxbasis{l}{i,j} |1\leq l\leq L_i^{(j)}\}.
\end{equation}
Furthermore, we define the local auxiliary space as
\begin{equation}
\aux^j = \aux^{(1,j)}\times\aux^{(2,j)}.
\end{equation}
The global auxiliary space is defined as the direct sum of all local auxiliary spaces as 
\begin{equation}
	\aux = \oplus_{j = 1}^N \aux^j.
\end{equation}

\textbf{Step 2. Construction of multiscale basis functions.} 
For the convenience of describing the method and the subsequent convergence analysis, for all $v = (v_1,v_2)\in \lsquare{K_j}$, we introduce a local projection operator $\pi_j:\lsquare{K_j}\rightarrow\aux$ as
\begin{equation}
\pi_j(v) = 
\left(
\sum_{l = 1}^{L_1^{(j)}}(v_1,\auxbasis{l}{1,j})\auxbasis{l}{1,j},
\sum_{l = 1}^{L_2^{(j)}}(v_2,\auxbasis{l}{2,j})\auxbasis{l}{2,j}
\right)
\end{equation}
and a global projection operator $\pi: \lsquare{\Omega}\rightarrow\aux$ is defined as 
\begin{equation}
\pi(v) = \sum_{j=1}^N\pi_j(v),\quad \forall v \in \lsquare{\Omega}.
\label{eq:basis_proj}
\end{equation}

The local multiscale basis functions are constructed by the following variational form of minimization problem. 
Within the oversampling domain $\osdomain{j}$ of every $K_j\subset \Omega$, find $\msbasis{l}{i,j} \in V(\osdomain{j})$ and $\overline{T}^{i,j,l}_{i',j',l'} \in \R$ such that for $i=1,2$, we have
\begin{equation}
	\begin{split}
	a_Q(\msbasis{l}{i,j},w)+ \sum_{i'=1}^2\sum_{\kij{l'}{i',j'}\subset\osdomain{j}}\overline{T}^{i,j,l}_{i',j',l'}(w\cdot \unit{i'},\auxbasis{l'}{i',j'}) & = 0, \quad \forall w\in V(\osdomain{j}),\\
	(\msbasis{l}{i,j}\cdot \unit{i'},\auxbasis{l'}{i',j'}) &= \delta_{i,i'}\delta_{l, l'}\delta_{j, j'},\quad \forall \kij{l'}{i',j'}\subset\osdomain{j}.
	\end{split}
	\label{eq:ms_min}
\end{equation}
Here, $\unit{i}$ is the canonical basis for $\R^2$. 
$\delta_{i,i'}$, $\delta_{l, l'}$ and $\delta_{j, j'}$ are the delta Dirac function. 
We use the local multiscale basis functions to obtain the multiscale finite element space, which will be used for deriving the multiscale solution, as 
\begin{equation}
	\msspace = \mathrm{span}\{\msbasis{l}{i,j} |1\leq l\leq L_i^{(j)}, 1\leq j \leq N, i = 1,2\}.
\end{equation}

The local multiscale basis functions are inspired from the global multiscale basis functions which are constructed in the similar way but on the global domain. For every coarse element $K_j$ in $\Omega$, find $\globasis{l}{i,j}\in V$ and $T^{i,j,l}_{i',j',l'}\in\R$ such that
\begin{equation}
	\begin{split}
	a_Q(\globasis{l}{i,j},w)+ \sum_{i'=1}^2\sum_{\kij{l'}{i',j'}\subset\Omega}T^{i,j,l}_{i',j',l'}(w\cdot \unit{i'},\auxbasis{l'}{i',j'}) & = 0, \quad \forall w\in V,\\
	(\globasis{l}{i,j}\cdot \unit{i'},\auxbasis{l'}{i',j'}) &= \delta_{i,i'}\delta_{l, l'}\delta_{j, j'},\quad \forall \kij{l'}{i',j'}\subset\Omega.
	\end{split}
	\label{eq:glo_min}
\end{equation}
The global multiscale finite element space is thus defined as
\begin{equation}
	\glo = \mathrm{span}\{\globasis{l}{i,j} | 1\leq l\leq \lij, 1\leq j \leq N, i = 1,2\}.
\end{equation}
As our analysis suggest, the global basis functions exhibit exponential decay properties and have small values outside a sufficiently large oversampling region. The fact suggests that we can save computational costs by localizing the basis functions by truncating the domain, without introducing huge error. 

We remark that if we denote $\tilde{V}$ as the null space of the global projection operator $\pi$, for any $\globasis{l}{i,j}\in \glo$, we have
\begin{equation}
	a_Q(\globasis{l}{i,j},v) = 0, \quad \forall v\in\tilde{V}.
\end{equation}
This implies that with respect to the inner product of $a_Q$, $\tilde{V}\subset \glo^{\bot}$. 
As a matter of fact, $\tilde{V} = \glo^{\bot}$.
 
\textbf{Step 3. Multiscale solution.} 
The process of finding the multiscale solution can be described as follows. 
Find $p_{\text{ms}} = (p_{\text{ms},1}, p_{\text{ms},2})$ with $p_{\text{ms}}(t,\cdot) \in \msspace$ s.t. for all $v = (v_1, v_2)$ with $v(t, \cdot) \in \msspace$,
\begin{equation}
c \left(\dfrac{\partial p_{\text{ms}}}{\partial t}, v \right) + a_Q(p_{\text{ms}}, v) = (f,v).
\label{eq:sol_ms}
\end{equation}

\section{Convergence Analysis}\label{sec:analysis}
In this section, we will analyze the proposed method.
First, we define the following norms and semi-norms on $V$:
\begin{equation}
\begin{split}
\norm{p}_c^2 & = c(p,p), \\ 
\norm{p}_a^2 & = a(p,p), \\
\vert p \vert_q^2 & = q(p,p), \\
\norm{p}_{a_Q}^2 & = a_Q(p,p),\\
\norm{p}_{\lweight{\Omega}{\kappa}}^2 & = \sum_i (\kappa_i^{\frac{1}{2}}p_i,\kappa_i^{\frac{1}{2}}p_i),\\
\norm{p}_{\lweight{\Omega}{\kappa^{-1}}}^2 & = \sum_i (\kappa_i^{-\frac{1}{2}}p_i,\kappa_i^{-\frac{1}{2}}p_i).
\end{split}
\end{equation}
For a subdomain $D = \bigcup_{j \in J} K_j$ as a union of coarse grid blocks, we also define the following local norms and semi-norms on $V$:
\begin{equation}
\begin{split}
\norm{p}_{a(D)}^2 & = \sum_{j \in J} a^{(j)}(p,p), \\
\vert p \vert_{q(D)}^2 & = \sum_{j \in J} q^{(j)}(p,p), \\
\norm{p}_{a_Q(D)}^2  & = \sum_{j \in J} a_Q^{(j)}(p,p),\\
\norm{p}_{\lweight{D}{\kappa}}^2 & = \sum_{j \in J} (\kappa_i^{\frac{1}{2}}p_i,\kappa_i^{\frac{1}{2}}p_i)_{L^2(K_j)},\\
\norm{p}_{\lweight{D}{\kappa^{-1}}}^2 & = \sum_{j \in J} (\kappa_i^{-\frac{1}{2}}p_i,\kappa_i^{-\frac{1}{2}}p_i)_{L^2(K_j)}.
\end{split}
\end{equation}

We remark that 
\begin{equation}
\begin{split}
	\norm{p}_{\lweight{D}{\kappa}} &\leq \overline{\kappa}^{\frac{1}{2}}\norm{p}_{\lsquare{D}},\\
	\norm{p}_{\lweight{D}{\kappa^{-1}}} &\leq \underline{\kappa}^{-\frac{1}{2}}\norm{p}_{\lsquare{D}}.
\end{split}
	\label{eq:w2s}
\end{equation}

In addition, we introduce some operators which will be used in our analysis, 
namely $\Rglo: V \to \glo$ given by: 
for any $u \in V$, the image $\Rglo u \in \glo$ is defined by
\begin{equation}
a_Q(\Rglo u , v) = a_Q(u,v), \quad \forall v \in V_{glo},
\label{eq:elliptic_proj}
\end{equation}
and similarly, $\Rms: V \to \msspace$ given by: 
for any $u \in V$, the image $\Rms u \in \msspace$ is defined by
\begin{equation}
a_Q(\Rms u , v) = a_Q(u,v), \quad \forall v \in \msspace.
\label{eq:elliptic_proj_ms}
\end{equation}
We also define $\mathcal{C}: V \to V$ given by:
for any $u \in V$, the image $\mathcal{C}u \in V$ is defined by
\begin{equation}
(\mathcal{C}u , v) = c(u,v),\quad \forall v \in V.
\end{equation}
Moreover, the operator $\mathcal{A}: D(\mathcal{A}) \to \lsquare{\Omega}$ is 
defined on a subspace $D(\mathcal{A}) \subset V$ by:
for any $u \in D(\mathcal{A})$, the image $\mathcal{A}u \in \lsquare{\Omega}$ is defined by
\begin{equation}
(\mathcal{A}u , v) = a_Q(u,v), \quad \forall v \in V.
\end{equation}

The following lemma shows that the projection operator $\Rglo$ has a good approximation property with respect to the $a_Q$-norm and $L^2$-norm.

\begin{lemma}
\label{lemma1}
	Let $u\in D(\mathcal{A})$, then we have $u-\Rglo u \in \tilde{V}$ and 
	\begin{equation}
		\norm{u-\Rglo u}_{a_Q} \leq CH\norm{\mathcal{A}u}_{\lweight{\Omega}{\kappa^{-1}}}.
		\label{eq:lemma1.0a}
	\end{equation}
	and
	\begin{equation}
		\norm{u-\Rglo u}_{\lsquare{\Omega}}\leq CH^2\underline{\kappa}^{-\frac{1}{2}}\norm{\mathcal{A}u}_{\lweight{\Omega}{\kappa^{-1}}}.
		\label{eq:lemma1.0b}
	\end{equation}

	\begin{proof}
		By \eqref{eq:elliptic_proj}, we directly get $u-\Rglo u \in \glo$. This yields
		\begin{equation}
			a_Q(u-\Rglo u, \Rglo u) = 0.
			\label{eq:lemma1.1}
		\end{equation}
		Thus, we have 
		\begin{equation}
			\begin{split}
			a_Q(u-\Rglo u, u-\Rglo u)
			& = a_Q(u-\Rglo u, u) - a_Q(u-\Rglo u, \Rglo u)\\
			& = a_Q(u-\Rglo u, u)\\
			& = a_Q(u,u-\Rglo u) \\
			& = (\mathcal{A}u, u-\Rglo u)\\
			& \leq \norm{\mathcal{A}u}_{\lweight{\Omega}{\kappa^{-1}}}
			\norm{(u-\Rglo u)}_{\lweight{\Omega}{\kappa}}.
			\end{split}
			\label{eq:lemma1.2}
		\end{equation}
		Since $u-\Rglo u \in \tilde{V}$, we have $\pi_j(u-\Rglo u) = 0$ for all $l = 1,2,\ldots, L_j$ and $j = 1,2, \ldots, N$. 
		Moreover, the Poincar\'{e} inequality gives
		\begin{equation}
			\int_{\kij{l}{i,j}}[(u-\Rglo u)\cdot\unit{i}]^2 \leq CH^2 \int_{\kij{l}{i,j}}|\nabla[(u - \Rglo u)\cdot\unit{i}]|^2.
			\label{eq:lemma1.3}
		\end{equation}
		This yields that
		\begin{equation}
			\begin{split}
			\norm{(u-\Rglo u)}_{\lweight{\Omega}{\kappa}}^2 
			& = \sum_i\norm{\kappa_i^{\frac{1}{2}}(u-\Rglo u)\cdot \unit{i}}^2_{L^2(\Omega)}\\
			& = \sum_i\sum_{\kij{l}{i,j}\subset\Omega}\norm{\kappa_i^{\frac{1}{2}}(u - \Rglo u)\cdot \unit{i}}_{L^2(\kij{l}{i,j})}^2\\	
			& \leq CH^2\norm{u-\Rglo u}_{a_Q}^2.
			\end{split}	
			\label{eq:lemma1.4}
		\end{equation}
		Thus, we have
		\begin{equation}
			\norm{u-\Rglo u}_{a_Q}^2\leq CH\norm{\mathcal{A}u}_{\lweight{\Omega}{\kappa^{-1}}}\norm{u-\Rglo u}_{a_Q},
			\label{eq:lemma1.5}
		\end{equation}
		which gives the estimate in the \eqref{eq:lemma1.0a}. 
		
		The proof of \eqref{eq:lemma1.0b} follows a duality argument. 
		Define $w\in V$ such that
		\begin{equation}
			a_Q(w,v) = (u-\Rglo u, v)\quad \forall v\in V.
			\label{eq:lemma1.6}
		\end{equation}
		Then we have
		\begin{equation}
			\norm{u-\Rglo u}_{\lsquare{\Omega}}^2 = (u-\Rglo u, u-\Rglo u) = a_Q(w, u-\Rglo u).
			\label{eq:lemma1.7}
		\end{equation}
		Taking $v = \Rglo w\in\glo$ in \eqref{eq:elliptic_proj}, we obtain
		\begin{equation}
			a_Q(u-\Rglo u, \Rglo w) = 0.
			\label{eq:lemma1.8}
		\end{equation}
		Since $w\in D(\mathcal{A})$ and $\mathcal{A}w = u-\Rglo u$, we have
		\begin{equation}
			\begin{split}
			\norm{u-\Rglo u}_{\lsquare{\Omega}}^2
			& = a_Q(w-\Rglo w, u-\Rglo u)\\
			& \leq \norm{w-\Rglo w}_{a_Q}\norm{u-\Rglo u}_{a_Q}\\
			& \leq \left(CH\norm{\mathcal{A}w}_{\lweight{\Omega}{\kappa^{-1}}}\right) \left(CH\norm{\mathcal{A}u}_{\lweight{\Omega}{\kappa^{-1}}}\right)\\
			& \leq \left(CH\underline{\kappa}^{-\frac{1}{2}}\norm{\mathcal{A}w}_{\lsquare{\Omega}}\right) \left(CH\norm{\mathcal{A}u}_{\lweight{\Omega}{\kappa^{-1}}}\right)\\
			& \leq CH^2\underline{\kappa}^{-\frac{1}{2}}\norm{u-\Rglo u}_{\lsquare{\Omega}}
						\norm{\mathcal{A}u}_{\lweight{\Omega}{\kappa^{-1}}}.
			\label{eq:lemma1.9}
			\end{split}
		\end{equation}
	\end{proof}
\end{lemma}

Next, we show that the global basis functions are localizable. 
For the purpose of this, for each coarse block $K$, we define a bubble function $B$ such that $B(x)>0, \forall x\in \interior(K)$ and $B(x) = 0, \forall x\in \partial K$. 
We will take $B = \prod_{x_k}\chi_k^H$, where $\chi_k^H$ is coarse scale partition of unity on $K$. 
Based on the bubble function, we define a constant as follows. 
\begin{equation}
C_{\text{equiv}} = \sup_{K_j\in\mesh^H,v\in\aux}\dfrac{\norm{v}_{\lsquare{K_j}}^2}{\norm{B^{\frac{1}{2}}v}_{\lsquare{K_j}}^2}.
\label{eq:cons_equiv}
\end{equation}

\begin{lemma}
\label{lemma2}
	For all $v_{\text{aux}}\in\aux$, there exists a function $v\in V$ such that
	\begin{equation}
		\pi(v) = v_{\text{aux}},
		\quad 
		\norm{v}_{a_Q}^2\leq D\norm{v_{\text{aux}}}_{\lweight{\Omega}{\kappa}}^2,
		\quad 
		\supp(v)\subset \supp(v_{\text{aux}}),
		\label{eq:lemma2.0}
	\end{equation}
	where $D = \frac{C_\mathcal{T}^2}{H^2}+2\max_i\norm{\frac{\sigma_i}{\kappa_i}}_{L^{\infty}(\Omega)}$ and $C_\mathcal{T}$ is the maximum of vertices over all coarse elements. 
	\begin{proof}
		Without loss of generality we assume $v_{\text{aux}}\in\aux^j$ with $\norm{v_{\text{aux}}}_{\lsquare{K_j}}=1$. 
		We consider the following saddle point problem:
		find $v\in V_0(K_j)$ and $T^{i'}_{l'}\in\R$ such that 
		\begin{equation}
			\begin{split}
			a_Q(v,w)+ 
			\sum_{i'=1}^2 \sum_{\kij{l'}{i',j}\subset K_j} T^{i'}_{l'} (w \cdot \unit{i'},\auxbasis{l'}{i',j}) & = 0, 
			\quad \forall w \in V_0(K_j),\\
			((v-v_{\text{aux}}) \cdot \unit{i'},\auxbasis{l'}{i',j}) &= 0,
			\quad \forall \kij{l'}{i',j}\subset K_j.
			\end{split}
			\label{eq:lemma2.2}
		\end{equation}
		The well-posedness of the above saddle point problem is equivalent to the existence of $\widetilde{v}\in V_0(K_j)$ such that
		\begin{equation}
			(\widetilde{v},v_{\text{aux}}) \geq C_1\norm{v_{\text{aux}}}_{\lweight{K_j}{\kappa}}^2,
			\quad \norm{\widetilde{v}}_{a_Q(K_j)} \leq C_2\norm{v_{\text{aux}}}_{\lweight{K_j}{\kappa}},
			\label{eq:lemma2.3}
		\end{equation}
		where $C_1$, $C_2$ are constants to be determined. 
		Taking $\widetilde{v} = B v_{\text{aux}}$, we have 
		\begin{equation}
			(\widetilde{v},v_{\text{aux}}) 
			= \norm{B^{\frac{1}{2}}v_{\text{aux}}}_{\lsquare{K_j}}^2
			\geq C_{\text{equiv}}^{-1}\norm{v_{\text{aux}}}_{\lsquare{K_j}}^2.
			\label{eq:lemma2.4}
		\end{equation}
		On the other hand, on every $\kij{l}{i,j}$, we have
		\begin{equation}
			\nabla(\widetilde{v}\cdot\unit{i}) = B\nabla(\widetilde{v}\cdot\unit{i})+\nabla B(\widetilde{v}\cdot\unit{i}).
			\label{eq:lemma2.5}
		\end{equation}
		By definition of $\aux$, $\norm{B\nabla(\widetilde{v}\cdot\unit{i})}_{L^2(K_j)}=0$. 
		At the same time, $|B|\leq 1$, $|\nabla B|\leq C_\mathcal{T} H^{-1}$. 
		Thus, 
		\begin{equation}
			\norm{\kappa_i^{\frac{1}{2}}\nabla(\widetilde{v}\cdot\unit{i})}_{L^2(K_j)}^2 \leq \dfrac{C_\mathcal{T}^2}{H^2}\norm{\kappa_i(\widetilde{v}\cdot\unit{i})}_{L^2(K_j)}^2.
			\label{eq:lemma2.6}
		\end{equation}
		This yields
		\begin{equation}
			\norm{\widetilde{v}}_{a(K_j)}^2 \leq 
			\dfrac{C_\mathcal{T}^2}{H^2}\sum_i \norm{\kappa_i^\frac{1}{2}(v_{\text{aux}}\cdot \unit{i})}_{L^2(K_j)}^2
			\label{eq:lemma2.7}
		\end{equation}
		and
		\begin{equation}
			|\widetilde{v}|_{q(K_j)}^2 \leq 
			2 \max_i \norm{\dfrac{\sigma_i}{\kappa_i}}_{L^{\infty}(K_j)}\sum_i\norm{\kappa_i^{\frac{1}{2}}(v_{\text{aux}}\cdot\unit{i})}_{L^2(K_j)}^2.
			\label{eq:lemma2.8}
		\end{equation}
		Thus, we have
		\begin{equation}
			\norm{\widetilde{v}}_{a_Q(K_j)}^2 \leq \left(\frac{C_\mathcal{T}^2}{H^2}+2\max_i\norm{\frac{\sigma_i}{\kappa_i}}_{L^{\infty}(\Omega)}\right) \norm{v_{\text{aux}}}_{\lweight{K_j}{\kappa}}^2,
			\label{eq:lemma2.9}
		\end{equation}
		This guarantees the existence and uniqueness of $v\in V_0(K_j)$ and $T^{i'}_{l'}\in\R$ satisfying \eqref{eq:lemma2.2}, in which $v$ satisfies our desired properties. 
	\end{proof}
\end{lemma}

Here, we make a remark that we can assume $D \geq 1$ without loss of generality.

In order to estimate the difference between the global basis functions and localized basis functions, 
we need the notion of a cutoff function with respect to the oversampling regions. 
For each coarse grid $K_j$ and $M > m$, we define $\chi_j^{M,m} \in \mathrm{span}\{ \chi_k^{H} \}$ 
such that $0 \leq \chi_j^{M,m} \leq 1$ and $\chi_j^{M,m} = 1$ on the inner region $K_{j,m}$ and 
$\chi_j^{M,m} = 0$ outside the region $K_{j,M}$.

The following lemma shows that our multiscale basis functions have a decay property. In particular, the
global basis functions are small outside an oversampling region specified in the lemma, 
which is important in localizing the multiscale basis functions.

\begin{lemma}
\label{lemma3}
Given $\auxbasis{l}{i,j} \in \aux^{j}$ and an oversampling region $K_{j,m}$ with number of layers $m \geq 2$. 
Let $\msbasis{l}{i,j}$ be a localized multiscale basis function defined on $\osdomain{j}$ given by \eqref{eq:ms_min}, and $\globasis{l}{i,j}$ be the corresponding global basis function given by \eqref{eq:glo_min}. 
Then we have
\begin{equation}
\norm{ \globasis{l}{i,j} - \msbasis{l}{i,j} }_{a_Q}^2 \leq E \norm{ \auxbasis{l}{i,j} }_{\lweight{K_j}{\kappa}}^2,
\end{equation}
where $E = 8D(2+D)(1+CH^2)\left(1+(C^\frac{1}{2} D^\frac{1}{2} H + CH^2)^{-1}\right)^{1-m}$.

\begin{proof}
By Lemma~\ref{lemma2}, there exists $v \in V$ such that
\begin{equation}
\pi(v) = \auxbasis{l}{i,j}, \quad 
\norm{ v }_{a_Q}^2 \leq D \norm{ \auxbasis{l}{i,j} }_{\lweight{\Omega}{\kappa}}^2, \quad 
\supp(v) \subset K_j.
\label{eq:lemma3.0}
\end{equation}
We take $\eta = \globasis{l}{i,j} - v \in V$
and $\zeta = v - \msbasis{l}{i,j} \in V(K_{j,m})$. 
Then $\pi(\eta) = \pi(\zeta) = 0$ and hence $\eta, \zeta \in \widetilde{V}$. 
We first see that for $K_{j'} \subset K_{j,m-1}$, 
\begin{equation}
\pi_{j'}(\chi_j^{m,m-1} \eta) = \pi_{j'}(\eta) = 0,
\end{equation}
since $\chi_j^{m,m-1} = 1$ on $K_{j,m-1}$ and $\eta \in \widetilde{V}$. 
On the other hand, for $K_{j'} \subset \Omega \setminus K_{j,m}$, 
\begin{equation}
\pi_{j'}(\chi_j^{m,m-1} \eta) = \pi_{j'}(0) = 0,
\end{equation}
since $\chi_j^{m,m-1} = 0$ on $\Omega \setminus K_{j,m}$. 
Therefore, we have $\text{supp}\left(\pi(\chi_j^{m,m-1} \eta)\right) \subset K_{j,m} \setminus K_{j,m-1}$. 
Again, by Lemma~\ref{lemma2}, there exists $\beta \in V$ such that
\begin{equation}
\pi(\beta) = \pi(\chi_j^{m,m-1} \eta), \quad 
\norm{ \beta }_{a_Q}^2 \leq D \norm{ \pi(\chi_j^{m,m-1} \eta) }_{\lweight{K_{j,m} \setminus K_{j,m-1}}{\kappa}}^2, \quad 
\supp(\beta) \subset K_{j,m} \setminus K_{j,m-1}.
\label{eq:lemma3.1}
\end{equation}
Take $\tau = \beta - \chi_j^{m,m-1} \eta \in V(K_{j,m})$. 
Again, $\pi(\tau) = 0$ and hence $\tau \in \widetilde{V}$.
Now, by the variational problems \eqref{eq:glo_min} and \eqref{eq:ms_min}, we have
\begin{equation}
\begin{split}
a_Q(\globasis{l}{i,j},w)+ \sum_{i'=1}^2\sum_{\kij{l'}{i',j'}\subset\Omega}T^{i,j,l}_{i',j',l'}(w\cdot \unit{i'},\auxbasis{l'}{i',j'}) & = 0, \quad \forall w\in V,\\
a_Q(\msbasis{l}{i,j},w)+ \sum_{i'=1}^2\sum_{\kij{l'}{i',j'}\subset\osdomain{j}}\overline{T}^{i,j,l}_{i',j',l'}(w\cdot \unit{i'},\auxbasis{l'}{i',j'}) & = 0, \quad \forall w\in V(\osdomain{j}) \\
\end{split}
\end{equation}
Taking $w = \tau - \zeta \in V(K_{j,m})$ and using the fact that $\tau - \zeta \in \widetilde{V}$, we have
\begin{equation}
a_Q(\globasis{l}{i,j} - \msbasis{l}{i,j}, \tau - \zeta) = 0,
\end{equation}
which implies
\begin{equation}
\begin{split}
\norm{\globasis{l}{i,j} - \msbasis{l}{i,j}}_{a_Q}^2 
& = a_Q(\globasis{l}{i,j} - \msbasis{l}{i,j}, \globasis{l}{i,j} - \msbasis{l}{i,j}) \\
& = a_Q(\globasis{l}{i,j} - \msbasis{l}{i,j}, \eta+\zeta) \\
& = a_Q(\globasis{l}{i,j} - \msbasis{l}{i,j}, \eta+\tau) \\
& \leq \norm{ \globasis{l}{i,j} - \msbasis{l}{i,j} }_{a_Q} \norm{ \eta+\tau }_{a_Q}.
\end{split}
\end{equation}
By \eqref{eq:lemma3.1}, we have
\begin{equation}
\begin{split}
\norm{\globasis{l}{i,j} - \msbasis{l}{i,j}}_{a_Q}^2 
& \leq \norm{ \eta+\tau}_{a_Q}^2 \\
& = \norm{ (1 - \chi_j^{m,m-1}) \eta + \beta }_{a_Q}^2 \\
& \leq 2\left(\norm{ (1 - \chi_j^{m,m-1}) \eta}_{a_Q}^2 + \norm{ \beta }_{a_Q}^2\right) \\
& \leq 2\left(\norm{ (1 - \chi_j^{m,m-1}) \eta}_{a_Q}^2 + D \norm{  \chi_j^{m,m-1} \eta }_{\lweight{K_{j,m} \setminus K_{j,m-1}}{\kappa}}^2 \right).
\end{split}
\label{eq:lemma3.2}
\end{equation}
For the first term on the right hand side of \eqref{eq:lemma3.2}, since 
\begin{equation}
\nabla \left((1 - \chi_j^{m,m-1}) (\eta \cdot \mathbf{e}_i) \right) = (1 - \chi_j^{m,m-1}) \nabla (\eta \cdot \unit{i}) - (\eta \cdot \unit{i}) \nabla \chi_j^{m,m-1},
\end{equation}
and $\vert 1 - \chi_j^{m,m-1} \vert \leq 1$, we have
\begin{equation}
\norm{ (1 - \chi_j^{m,m-1}) \eta }_{a}^2
\leq 2 \left( \norm{ \eta }^2_{a(\Omega \setminus K_{j,m-1})} + \norm{ \eta }_{\lweight{\Omega \setminus K_{j,m-1}}{\kappa}}^2 \right).
\end{equation}
On the other hand, we have
\begin{equation}
\vert (1 - \chi_j^{m,m-1}) \eta \vert_q^2 \leq \vert \eta \vert_{q(\Omega \setminus  K_{j,m-1})}^2.
\end{equation}
Therefore, we arrive at 
\begin{equation}
\norm{ (1 - \chi_j^{m,m-1}) \eta }_{a_Q}^2
\leq 2 \left( \norm{ \eta }^2_{a_Q(\Omega \setminus K_{j,m-1})} + \norm{ \eta }_{\lweight{\Omega \setminus K_{j,m-1}}{\kappa}}^2 \right).
\end{equation}
For the second term on the right hand side of \eqref{eq:lemma3.2}, 
using the fact that $\vert \chi_j^{m,m-1} \vert \leq 1$, we have
\begin{equation}
\begin{split}
\norm{ \pi(\chi_j^{m,m-1} \eta) }_{\lweight{K_{j,m} \setminus K_{j,m-1}}{\kappa}}^2 
& \leq \norm{ \chi_j^{m,m-1} \eta }_{\lweight{K_{j,m} \setminus K_{j,m-1}}{\kappa}}^2 \\
& \leq \norm{  \eta }_{\lweight{K_{j,m} \setminus K_{j,m-1}}{\kappa}}^2.
\end{split}
\end{equation}
To sum up, we have 
\begin{equation}
\norm{\globasis{l}{i,j} - \msbasis{l}{i,j}}_{a_Q}^2 
\leq 4 \norm{ \eta }^2_{a_Q(\Omega \setminus K_{j,m-1})} + (4+2D)  \norm{\eta }_{\lweight{\Omega \setminus K_{j,m-1}}{\kappa}}^2 .
\end{equation}
Since $\eta \in \widetilde{V}$, using Poincar\'{e} inequality, we obtain
\begin{equation}
\norm{ \eta }_{\lweight{\Omega \setminus K_{j,m-1}}{\kappa}}^2 \leq CH^2 \norm{\eta }^2_{a_Q(\Omega \setminus K_{j,m-1})}.
\end{equation}
Combining all the estimates, we have 
\begin{equation}
\norm{\globasis{l}{i,j} - \msbasis{l}{i,j}}_{a_Q}^2 
\leq (4+2D)(1+CH^2) \norm{ \eta }^2_{a_Q(\Omega \setminus K_{j,m-1})}.
\label{eq:lemma3.3}
\end{equation}

Next, we will prove a recursive estimate for $\norm{\eta}^2_{a_Q(\Omega \setminus K_{j,m-1})}$. 
We take $\xi = 1 - \chi_j^{m-1,m-2}$. 
Then $\xi = 1$ in $\Omega \setminus K_{j,m-1}$ and $0 \leq \xi \leq 1$. 
Hence, using 
\begin{equation}
\nabla(\xi^2 (\eta \cdot \unit{i})) = \xi^2 \nabla (\eta \cdot \unit{i})+ 2 \xi  (\eta \cdot \unit{i})\nabla \xi,
\end{equation}
we have
{
\begin{equation}
\vert \xi \eta \vert_a^2 
= a(\eta, \xi^2 \eta) + \norm{ \eta }_{\lweight{K_{j,m-1} \setminus K_{j,m-2}}{\kappa}}^2,
\end{equation}
}
which results in 
\begin{equation}
\norm{ \eta }_{a_Q(\Omega \setminus K_{j,m-1})}^2 
\leq \norm{ \xi \eta }_{a_Q}^2 
\leq a_Q(\eta, \xi^2 \eta)
+ {\| \eta \|_{\lweight{K_{j,m-1} \setminus K_{j,m-2}}{\kappa}}^2}.
\label{eq:lemma3.4}
\end{equation}
We will estimate the first term on the right hand side of \eqref{eq:lemma3.4}. 
Following the preceding argument, we see that $\supp(\pi(\xi^2 \eta)) \subset  K_{j,m-1} \setminus K_{j,m-2}$.
By Lemma~\ref{lemma2}, there exists $\gamma \in V$ such that
\begin{equation}
\pi(\gamma) = \pi(\xi^2 \eta), \quad 
\norm{ \gamma }_{a_Q}^2 \leq D \norm{ \pi(\xi^2 \eta) }_{\lweight{K_{j,m-1} \setminus K_{j,m-2}}{\kappa}}^2, \quad 
\supp(\gamma) \subset K_{j,m-1} \setminus K_{j,m-2}.
\label{eq:lemma3.5}
\end{equation}
Take $\theta = \xi^2 \eta - \gamma$. Again, $\pi(\theta) = 0$ and hence $\theta \in \widetilde{V}$. 
Therefore, we have 
\begin{equation}
a_Q(\globasis{l}{i,j}, \theta) = 0.
\end{equation}
Additionally, $\supp(\theta) \subset \Omega \setminus K_{j,m-2}$. 
Recall that, in \eqref{eq:lemma3.0}, we have $\supp(v) \subset K_j$.
Hence $\theta$ and $v$ have disjoint supports, and
\begin{equation}
a_Q(v, \theta) = 0.
\end{equation}
Therefore, we obtain
\begin{equation}
a_Q(\eta, \theta) = a_Q(\globasis{l}{i,j}, \theta) - a_Q(v, \theta) = 0.
\end{equation}
Note that $\xi^2 \eta = \theta + \gamma$. Using \eqref{eq:lemma3.5}, we have
\begin{equation}
\begin{split}
a_Q(\eta, \xi^2 \eta)
& = a_Q(\eta,\gamma) \\
& \leq \norm{ \eta }_{a_Q(K_{j,m-1} \setminus K_{j,m-2})} \norm{ \gamma }_{a_Q(K_{j,m-1} \setminus K_{j,m-2})} \\
& \leq D^{\frac{1}{2}} \norm{ \eta }_{a_Q(K_{j,m-1} \setminus K_{j,m-2})}  \norm{ \pi(\xi^2 \eta) }_{\lweight{K_{j,m-1} \setminus K_{j,m-2}}{\kappa}}.
\end{split}
\end{equation}
Since $\vert \xi \vert \leq 1$, we have 
\begin{equation}
\begin{split}
\norm{ \pi(\xi^2 \eta) }_{\lweight{K_{j,m-1} \setminus K_{j,m-2}}{\kappa}}
\leq \norm{ \xi^2 \eta }_{\lweight{K_{j,m-1} \setminus K_{j,m-2}}{\kappa}}
\leq \norm{ \eta }_{\lweight{K_{j,m-1} \setminus K_{j,m-2}}{\kappa}}.
\end{split}
\end{equation}
Hence, the right hand side of \eqref{eq:lemma3.4} can be estimated by 
\begin{equation}
\norm{ \eta }_{a_Q(\Omega \setminus K_{j,m-1})}^2 
\leq D^{\frac{1}{2}} \norm{ \eta }_{a_Q(K_{j,m-1} \setminus K_{j,m-2})}  \norm{ \eta }_{\lweight{K_{j,m-1} \setminus K_{j,m-2}}{\kappa}}
+ {\norm{ \eta}_{\lweight{K_{j,m-1} \setminus K_{j,m-2}}{\kappa}}^2}.
\label{eq:lemma3.6a}
\end{equation}
Since $\pi(\eta) = 0$, using Poincar\'{e} inequality, we have
\begin{equation}
\begin{split}
\norm{ \eta }^2_{\lweight{K_{j,m-1} \setminus K_{j,m-2}}{\kappa}}.
& \leq CH^2 \norm{ \eta }^2_{a_Q(K_{j,m-1} \setminus K_{j,m-2})},
\end{split}
\end{equation}
which implies
\begin{equation}
\norm{ \eta }^2_{a_Q(\Omega \setminus K_{j,m-1})} 
\leq (C^\frac{1}{2} D^\frac{1}{2} H + CH^2) \norm{\eta }^2_{a_Q(K_{j,m-1} \setminus K_{j,m-2})}.
\end{equation}
Therefore,
\begin{equation}
\begin{split}
\norm{ \eta }^2_{a_Q(\Omega \setminus K_{j,m-2})} 
& = \norm{ \eta }^2_{a_Q(\Omega \setminus K_{j,m-1})} + \norm{ \eta }^2_{a_Q( K_{j,m-1} \setminus K_{j,m-2})} \\
& \geq \left(1 + (C^\frac{1}{2} D^\frac{1}{2} H + CH^2)^{-1}\right) \norm{ \eta }^2_{a_Q(\Omega \setminus K_{j,m-1})}.
\end{split}
\end{equation}
Inductively, we have
\begin{equation}
\begin{split}
\norm{ \eta }^2_{a_Q(\Omega \setminus K_{j,m-1})} 
& \leq \left(1 +(C^\frac{1}{2} D^\frac{1}{2} H + CH^2)^{-1}\right)^{1-m} \norm{ \eta }^2_{a_Q(\Omega \setminus K_{j})} \\ 
& \leq \left(1 + (C^\frac{1}{2} D^\frac{1}{2} H + CH^2)^{-1}\right)^{1-m} \norm{ \eta }^2_{a_Q}.
\end{split}
\label{eq:lemma3.7}
\end{equation}
Finally, we estimate the term on the right hand side of \eqref{eq:lemma3.7}. 
Recall from the first property of $v$ in \eqref{eq:lemma3.0}, 
we have $\pi(v) = \auxbasis{l}{i,j}$, which implies 
\begin{equation}
(v \cdot \unit{i'},\auxbasis{l'}{i',j'}) = \delta_{i,i'}\delta_{l, l'}\delta_{j, j'},\quad \forall \kij{l'}{i',j'}\subset\Omega.
\end{equation}
Taking $w = \eta$ in \eqref{eq:glo_min}, we have 
\begin{equation}
a_Q(\globasis{l}{i,j}, \eta) = 0, 
\end{equation}
which implies 
\begin{equation}
\| \globasis{l}{i,j} \|_{a_Q} \leq \| v \|_{a_Q}.
\end{equation}
Using a triangle inequality and the second property of $v$ in \eqref{eq:lemma3.0}, we have 
\begin{equation}
\norm{ \eta }_{a_Q} = \norm{ \globasis{l}{i,j} - v }_{a_Q} \leq 2 \norm{ v }_{a_Q} \leq 2 D^\frac{1}{2}  \norm{ \auxbasis{l}{i,j} }_{\lweight{K_j}{\kappa}}.
\label{eq:lemma3.8}
\end{equation}
Combining \eqref{eq:lemma3.3}, \eqref{eq:lemma3.7} and \eqref{eq:lemma3.8}, we obtain our desired result.
\end{proof}
\end{lemma}

The following lemma shows that, similar to the global projection operator $\Rglo$, our localized multiscale finite element projection operator $\Rms$ can also provide a good approximation with respect to the $a_Q$-norm and $L^2$-norm.

\begin{lemma}
	\label{lemma4}
	Let $u\in D(\mathcal{A})$. Let $m \geq 2$ be the number of coarse grid layers in the oversampling regions in \eqref{eq:ms_min}. 
	If $m = O\left(\log\left(\frac{\overline{\kappa}}{H}\right)\right)$, we have
		\begin{equation}
			\norm{u-\Rms u}_{a_Q} \leq CH\norm{\mathcal{A}u}_{\lweight{\Omega}{\kappa^{-1}}}.
			\label{eq:lemma4.0a}
		\end{equation}
		and
		\begin{equation}
			\norm{u-\Rms u}_{\lsquare{\Omega}}\leq CH^2\underline{\kappa}^{-\frac{1}{2}}\norm{\mathcal{A}u}_{\lweight{\Omega}{\kappa^{-1}}}.
			\label{eq:lemma4.0b}
		\end{equation}
		
	\begin{proof}
		We write
		\begin{equation}
		\Rglo u = \sum_{i = 1}^2\sum_{j=1}^N\sum_{l=1}^{L_i^{(j)}}\alpha_l^{(i,j)}\globasis{l}{i,j}\in\glo
		\label{eq:lemma4.1}
		\end{equation}
		and define
		\begin{equation}
		w = \sum_{i = 1}^2\sum_{j=1}^N\sum_{l=1}^{L_i^{(j)}}\alpha_l^{(i,j)}\msbasis{l}{i,j}\in\msspace. 
		\label{eq:lemma4.2}
		\end{equation}
		By \eqref{eq:elliptic_proj_ms}, we have
		\begin{equation}
		\norm{u-\Rms u}_{a_Q}\leq\norm{u-w}_{a_Q}\leq \norm{u-\Rglo u}_{a_Q}+\norm{\Rglo u-w}_{a_Q}
		\label{eq:lemma4.3}
		\end{equation}
		By Lemma \ref{lemma3}, we have that
		\begin{equation}
		\begin{split}
		\norm{\Rglo u -w}_{a_Q}^2 
		& = \norm{\sum_{i = 1}^2\sum_{j=1}^N\sum_{l=1}^{L_i^{(j)}}\alpha_l^{(i,j)}(\globasis{l}{i,j}-\msbasis{l}{i,j})}_{a_Q}^2\\
		& \leq C(m+1)^2\sum_{j = 1}^N\norm{\sum_{i=1}^2\sum_{l=1}^{L_i^{(j)}}\alpha_l^{(i,j)}(\globasis{l}{i,j}-\msbasis{l}{i,j})}_{a_Q}^2\\
		& \leq CE(m+1)^2\sum_{j=1}^N\norm{\sum_{i=1}^2\sum_{l=1}^{L_i^{(j)}} \alpha_l^{(i,j)} \auxbasis{l}{i,j}}_{\lweight{K_j}{\kappa}}^2\\
		& \leq CE(m+1)^2\norm{\Rglo u}_{\lweight{\Omega}{\kappa}}^2.
		\end{split}
		\label{eq:lemma4.4}
		\end{equation}
		Combining \eqref{eq:lemma4.3}, \eqref{eq:lemma4.4} and Lemma \ref{lemma1}, we have
		\begin{equation}
		\norm{u-\Rms u}_{a_Q} \leq CH\norm{\mathcal{A}u}_{\lweight{\Omega}{\kappa^{-1}}} + CE^{\frac{1}{2}}(m+1)\norm{\Rglo u}_{\lweight{\Omega}{\kappa}}.
		\label{eq:lemma4.5}
		\end{equation}
		Now, we estimate $\norm{\Rglo u}_{\lweight{\Omega}{\kappa}}^2$. 
		By Ponicar\'{e} inequality, we have
		\begin{equation}
		\begin{split}
		\norm{\Rglo u}_{\lweight{\Omega}{\kappa}}^2
		\leq \overline{\kappa}\norm{\Rglo u}_{\lsquare{\Omega}}^2
		\leq C_p \overline{\kappa}\underline{\kappa}^{-1}\norm{\Rglo u}_{a_Q}^2.
		\end{split}
		\label{eq:lemma4.6}
		\end{equation}
		Taking $v = \Rglo u$ in \eqref{eq:elliptic_proj} and by Cauchy-Schwarz inequality, we obtain
		\begin{equation}
		\begin{split}
		\norm{\Rglo u}_{a_Q}^2 
		= a_Q(u,\Rglo u) 
		= (\mathcal{A}u,\Rglo u)
		\leq \norm{\mathcal{A}u}_{\lweight{\Omega}{\kappa^{-1}}}\norm{\Rglo u}_{\lweight{\Omega}{\kappa}}.
		\end{split}
		\label{eq:lemma4.7}
		\end{equation}
		Combining \eqref{eq:lemma4.6} and \eqref{eq:lemma4.7}, we obtain
		\begin{equation}
		\norm{\Rglo u}_{\lweight{\Omega}{\kappa}}\leq C\overline{\kappa}\underline{\kappa}^{-1}\norm{\mathcal{A}u}_{\lweight{\Omega}{\kappa^{-1}}}.
		\label{eq:lemma4.8}
		\end{equation}
		This yields
		\begin{equation}
		\norm{u-\Rms u }_{a_Q}\leq C(H+\overline{\kappa}\underline{\kappa}^{-1}E^{\frac{1}{2}}(m+1))\norm{\mathcal{A}u}_{\lweight{\Omega}{\kappa^{-1}}}
		\label{eq:lemma4.9}
		\end{equation}
		To obtain desired result, we will need
		\begin{equation}
		\label{eq:lemma4.10}
		H^{-1}\overline{\kappa}\underline{\kappa}^{-1}E^{\frac{1}{2}}(m+1) = O(1).
		\end{equation}
		Taking logarithm, we have
		\begin{equation}
		\log(H^{-1}) + \log(\overline{\kappa}) - \log(\underline{\kappa})+\frac{1-m}{2}\log\left(1+(C^\frac{1}{2} D^\frac{1}{2} H + CH^2)^{-1}\right) = O(1).
		\label{eq:lemma4.11}
		\end{equation}
		Therefore, if we take $m = O\left(\log\left(\frac{\overline{\kappa}}{H}\right)\right)$, we have \eqref{eq:lemma4.0a}. 
		The proof of \eqref{eq:lemma4.0b} follows a similar duality argument as in Lemma \ref{lemma1}. 
	\end{proof}
\end{lemma}

Now we are ready to obtain our main theorem on estimating error between $p$ and $\mssol$. 
\begin{theorem}
	\label{theorem}
	Suppose $f\in\lsquare{\Omega}$. 
	Let $m\geq 2$ be the number of coarse grid layers of the oversampling domain in \eqref{eq:ms_min}. 
	Let $p$ be the solution of \eqref{eq:sol_weak} and $\mssol$ of \eqref{eq:sol_ms}. 
	If $m = O\left(\log\left(\frac{\overline{\kappa}}{H}\right)\right)$, we have
	\begin{equation}
		\norm{p(T,\cdot)-\mssol(T,\cdot)}_c^2+\int_0^T\norm{p-\mssol}_{a_Q}^2\dif t 
		\leq CH^2\underline{\kappa}^{-1}\left(\norm{p^0}_{a_Q}^2+\int_0^T\norm{f}_{\lsquare{\Omega}}^2\dif t\right).
		\label{eq:theorem.0}
	\end{equation}
	
	\begin{proof}
		Taking $v = \frac{\partial p}{\partial t}$ in \eqref{eq:sol_weak}, we have
		\begin{equation}
			\norm{\frac{\partial p}{\partial t}}_c^2 
			+ \frac{1}{2}\frac{\dif}{\dif t}\norm{p}_{a_Q}^2 
			= \left(f,\frac{\partial p}{\partial t}\right)
			\leq C\norm{f}_{\lsquare{\Omega}}^2 + \frac{1}{2}\norm{\frac{\partial p}{\partial t}}_c^2.
		\label{eq:theorem.1}
		\end{equation}
		Integrating over $(0,T)$, we have
		\begin{equation}
			\frac{1}{2}\int_0^T\norm{\frac{\partial p}{\partial t}}_c^2 \dif t
			+ \frac{1}{2}\norm{p(T,\cdot)}_{a_Q}^2
			\leq C\left(\norm{p^0}_{a_Q}^2 + \int_0^T\norm{f}_{\lsquare{\Omega}}^2\dif t\right)
			\label{eq:theorem.2}
		\end{equation}
		Similarly, by taking $v = \frac{\partial \mssol}{\partial t}$ in \eqref{eq:sol_ms} and integrating over $(0,T)$, we obtain
		\begin{equation}
			\frac{1}{2}\int_0^T\norm{\frac{\partial \mssol}{\partial t}}_c^2 \dif t
			+ \frac{1}{2}\norm{\mssol(T,\cdot)}_{a_Q}^2
			\leq C\left(\norm{p^0}_{a_Q}^2 + \int_0^T\norm{f}_{\lsquare{\Omega}}^2\dif t\right)
			\label{eq:theorem.3}
		\end{equation}
		At the same time, by \eqref{eq:sol_weak}, we can see that
		\begin{equation}
			\mathcal{A}p = f - \mathcal{C}\frac{\partial p}{\partial t}
			\label{eq:theorem.4}.
		\end{equation}
		Thus, we have
		\begin{equation}
			\norm{\mathcal{A}p}_{\lsquare{\Omega}}^2
			\leq C\left(\norm{f}_{\lsquare{\Omega}}^2 + \norm{\frac{\partial p}{\partial t}}_c\right)
			\label{eq:theorem.5}
		\end{equation}
		By the definition of $p$ and $\mssol$ in \eqref{eq:sol_weak} and \eqref{eq:sol_ms}, respectively, we can get that $\forall v\in\msspace$ and $t\in(0,T)$, we have
		\begin{equation}
			c\left(\frac{\partial (p-\mssol)}{\partial t}, v\right) + a_Q(p-\mssol,v) = 0.
			\label{eq:theorem.6}
		\end{equation}
		Thus, we have
		\begin{equation}
			\begin{split}
			& \frac{1}{2}\frac{\dif}{\dif t}\norm{p-\mssol}_c^2+\norm{p-\mssol}_{a_Q}^2\\
			& = c\left(\frac{\partial(p-\mssol)}{\partial t},p-\mssol\right)+a_Q(p-\mssol,p-\mssol)\\
			& = c\left(\frac{\partial (p-\mssol)}{\partial t}, p-\Rms p\right)+ a_Q(p-\mssol,p-\Rms p)\\
			& \leq \norm{\frac{\partial (p-\mssol)}{\partial t}}_c\norm{p-\Rms p}_c
				 + \norm{p-\mssol}_{a_Q}\norm{p-\Rms p}_{a_Q}\\
			& \leq \left(\norm{\frac{\partial p}{\partial t}}_c + \norm{\frac{\partial \mssol}{\partial t}}_c\right)\norm{p-\Rms p}_c 
				+ \frac{1}{2}\norm{p-\mssol}_{a_Q}^2 + \frac{1}{2}\norm{p-\Rms p}_{a_Q}^2.
			\end{split}
			\label{eq:theorem.7}
		\end{equation}
		Integrating over $(0,T)$ and using \eqref{eq:theorem.4} by Lemma \ref{lemma4} with \eqref{eq:w2s}, we obtain
		\begin{equation}
			\begin{split}
			& \frac{1}{2}\norm{p(T,\cdot)-\mssol(T,\cdot)}_c^2 + \frac{1}{2}\int_0^T\norm{p-\mssol}_{a_Q}^2\dif t\\
			& \leq \int_0^T\left(\norm{\frac{\partial p}{\partial t}}_c + \norm{\frac{\partial \mssol}{\partial t}}_c\right)\norm{p-\Rms p}_c\dif t
				+ \frac{1}{2}\int_0^T\norm{p-\Rms p}_{a_Q}^2\dif t\\
			& \leq \left(\int_0^T\left(\norm{\frac{\partial p}{\partial t}}_c + \norm{\frac{\partial \mssol}{\partial t}}_c\right)^2\dif t\right)^{\frac{1}{2}} \left(\int_0^T\norm{p-\Rms p}_c^2\dif t\right)^{\frac{1}{2}} 
				+ \frac{1}{2}\int_0^T\norm{p-\Rms p}_{a_Q}^2\dif t\\
			& \leq \left(\int_0^T\left(\norm{\frac{\partial p}{\partial t}}_c + \norm{\frac{\partial \mssol}{\partial t}}_c\right)^2\dif t\right)^{\frac{1}{2}} \left(\int_0^T CH^4\underline{\kappa}^{-2}\left(\norm{f}_{\lsquare{\Omega}}+\norm{\frac{\partial p}{\partial t}}_c\right)^2\dif t\right)^{\frac{1}{2}}\\
				&\quad + \int_0^T CH^2\underline{\kappa}^{-1}\left(\norm{f}_{\lsquare{\Omega}}+\norm{\frac{\partial p}{\partial t}}_c\right)^2\dif t\\
			& \leq CH^2\underline{\kappa}^{-1}
				\int_0^T \left(\norm{\frac{\partial p}{\partial t}}_c^2 + \norm{\frac{\partial p}{\partial t}}_c^2 + \norm{f}^2_{\lsquare{\Omega}}\right) \dif t.
			\end{split}
			\label{eq:theorem.8}
		\end{equation}
		Combining \eqref{eq:theorem.2}, \eqref{eq:theorem.3} and \eqref{eq:theorem.8}, we get the result in the theorem.
	\end{proof}
\end{theorem}

\section{Numerical Examples}\label{sec:numerical}

In this section, we present two numerical examples with high-contrast media to 
verify the convergence of our proposed method,  
using a fine-scale approximation $p_{\text{f}}$ as a reference solution. 
We will compute the coarse cell average $\meansol_{\text{f}}$ of 
the fine-scale solution $p_{\text{f}}$ and $\meansol_{\text{ms}}$ of the multiscale solution $p_{\text{ms}}$, and
compare the relative $L_2$ error of coarse cell average, i.e. 
\begin{equation}\label{eq:error}
	e_{L_2}^{(i)} = \|\meansol_{\text{f},i}-\meansol_{\text{ms},i}\|_{L^2},
	\qquad \|\meansol_{\text{f},i}-\meansol_{\text{ms},i}\|^2_{L_2} = \dfrac{\sum_K(\meansol^K_{\text{f},i}-\meansol_{\text{ms},i}^K)^2}{\sum_K(\meansol^K_f)^2},
	\qquad \meansol_{\text{f},i}^K = \dfrac{1}{|K|}\int_K p_{\text{f},i} \dif x. 
\end{equation}
In all the experiments, we take the spatial domain to be $\Omega = (0,1)^2$ and 
the fine mesh size to be $h = 1/256$. 
An example of the media $\kappa_1$ and $\kappa_2$ used in the experiments is illustrated in Figure~\ref{fig:kappa1}. 
In the figure, the contrast values, i.e. the ratio of the maximum and the minimum in $\Omega$, of the media are
$\overline{\kappa}_1 = \overline{\kappa}_2 = 10^4$. 
Unless otherwise specified, we set $\sigma = 1$.

\begin{figure}[!ht]
	\centering
	\begin{subfigure}{.45\textwidth}
		\centering
		\includegraphics[width=.9\linewidth]{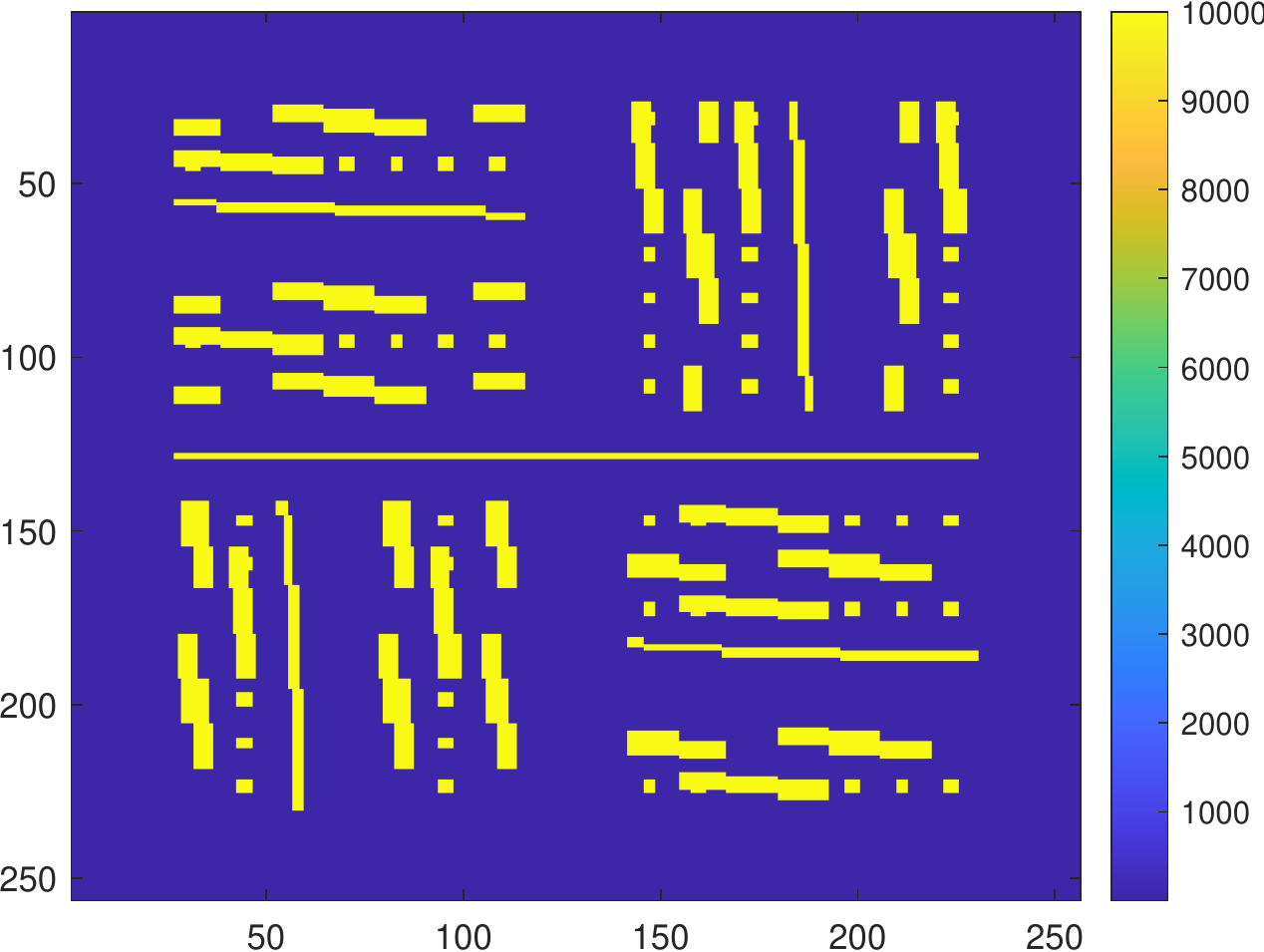}
		\caption{$\kappa_1$}
	\end{subfigure}
	\begin{subfigure}{.45\textwidth}
		\centering
		\includegraphics[width = .9\linewidth]{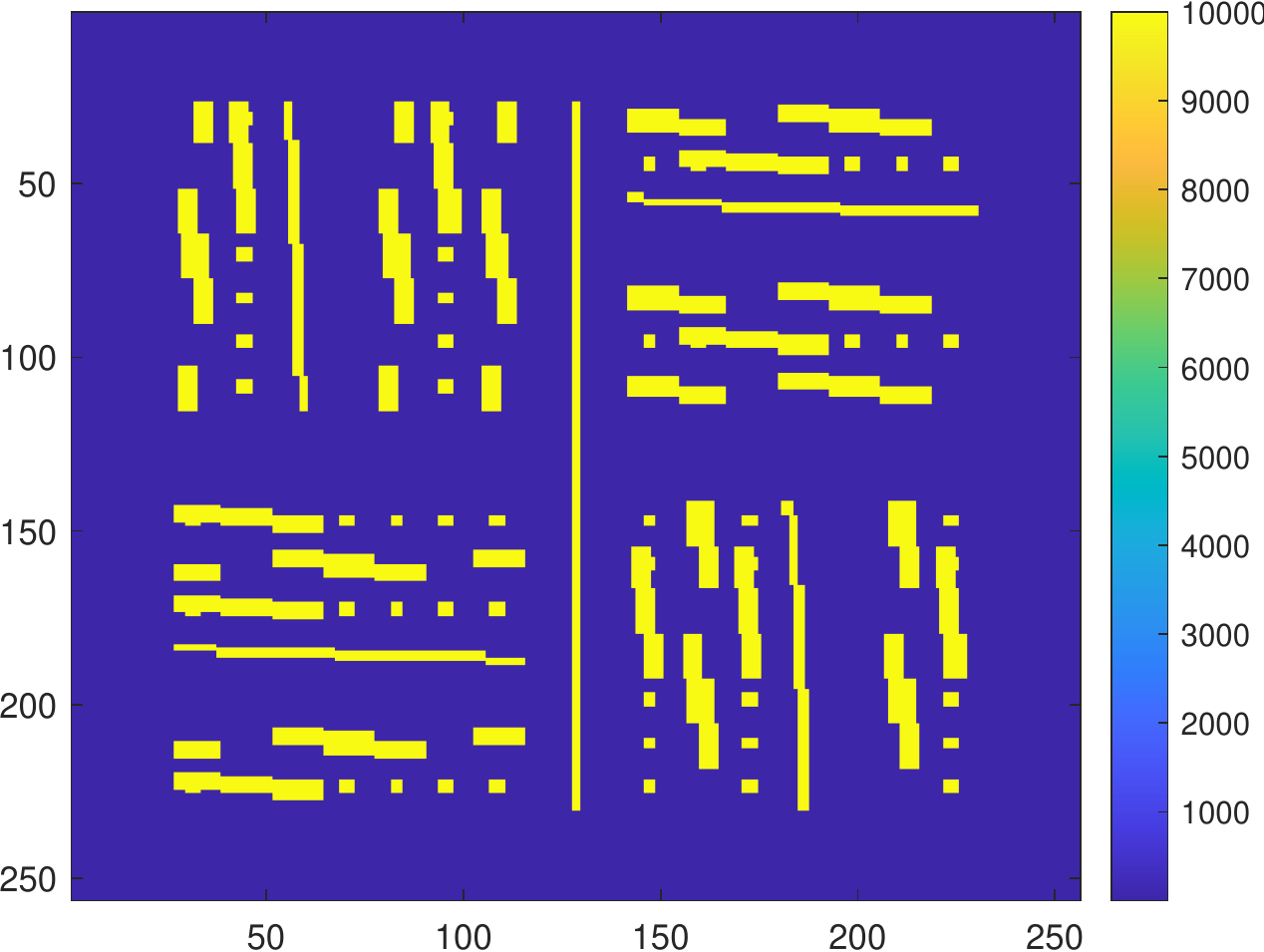}
		\caption{$\kappa_2$}
	\end{subfigure}
\caption{High contrast permeability field}
\label{fig:kappa1}
\end{figure}

\subsection{Experiment 1 static case}
In this experiment, we consider the dual continuum model under steady state
\begin{equation}
\begin{split}
 - \text{div}(\kappa_1 \nabla p_1) + \sigma(p_1 - p_2) = f_1, \\
 - \text{div}(\kappa_2 \nabla p_2) - \sigma(p_1 - p_2) = f_2. 
\end{split}
\label{eq:ss}
\end{equation} 
The source terms are given as $f_1(x,y) = 1$ and $f_2(x,y)$ as shown in Figure~\ref{fig:f2_static}. 
In Figure~\ref{fig:sol_static}, we plot the fine-scale solution, the coarse-scale average 
and the NLMC coarse-scale solution with coarse mesh size $H = 1/64$ and number of oversampling layers $m = 8$, 
from which we observe very good agreement between the coarse-scale average and the NLMC solution.
In Table~\ref{tab:conv_static}, we present the relative $L^2$ error with varying coarse grid size. 
With the number of oversampling layers satisfying the sufficient condition, 
we can see that the error converges. 
We also compare the performance of different numbers of oversampling layers under fixed coarse mesh size $H$. 
The results are summarized in Table~\ref{tab:error_32} for $H = 1/32$ and Table~\ref{tab:error_64} for $H = 1/64$. 
It can be seen that the error decays quickly with respect to the number of oversampling layers $m$ for both cases, 
which verifies the fact that the oversampling region has to be sufficiently large to obtain quality numerical approximations. 
\begin{figure}[!ht]
	\centering
	\includegraphics[width = .45\linewidth]{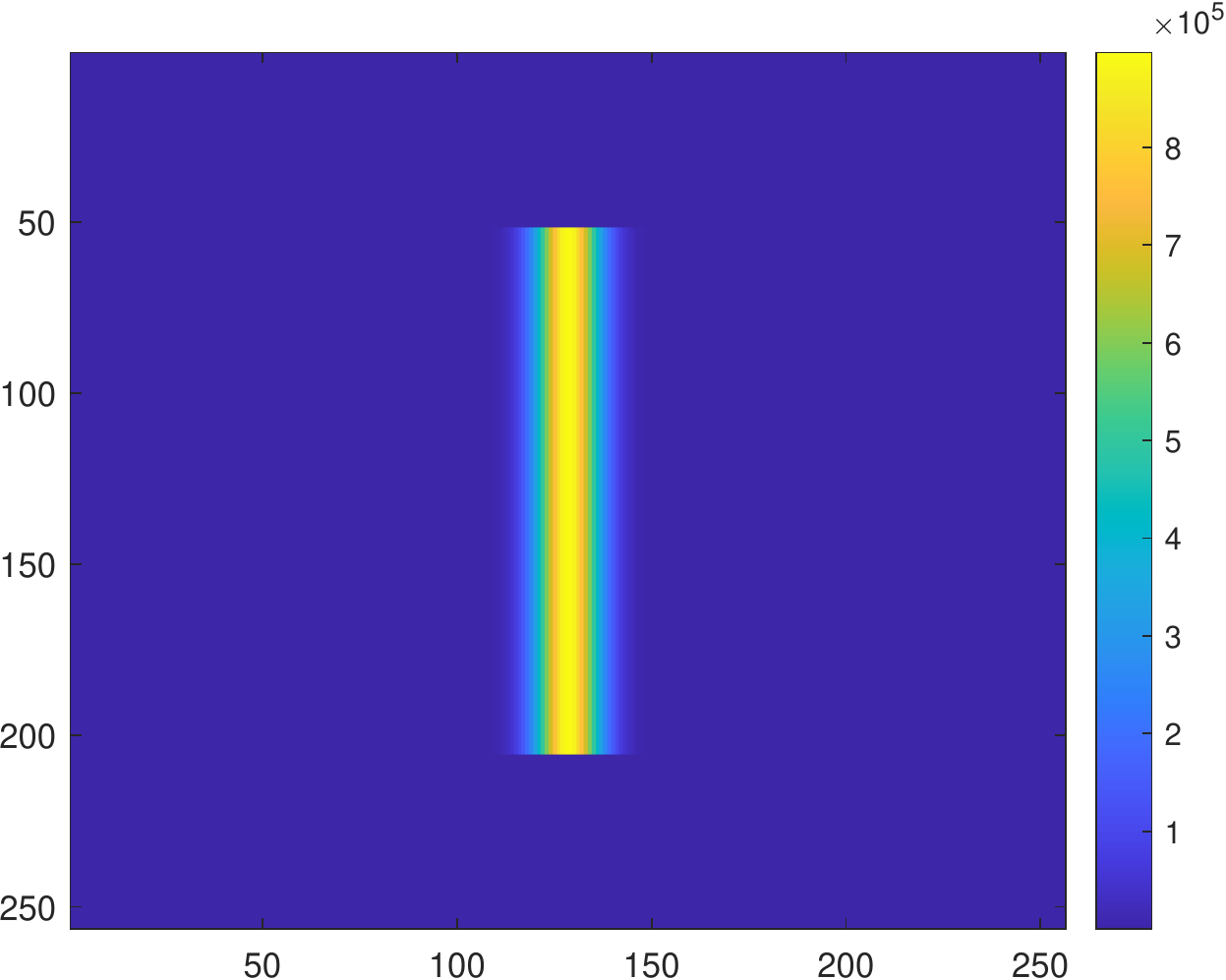}
	\caption{Source term $f_2$ in Experiment 1.}
	\label{fig:f2_static}
\end{figure}

%
%

\begin{figure}[!ht]
	\centering
	\begin{subfigure}{.45\textwidth}
		\centering
		\includegraphics[width = .9\linewidth]{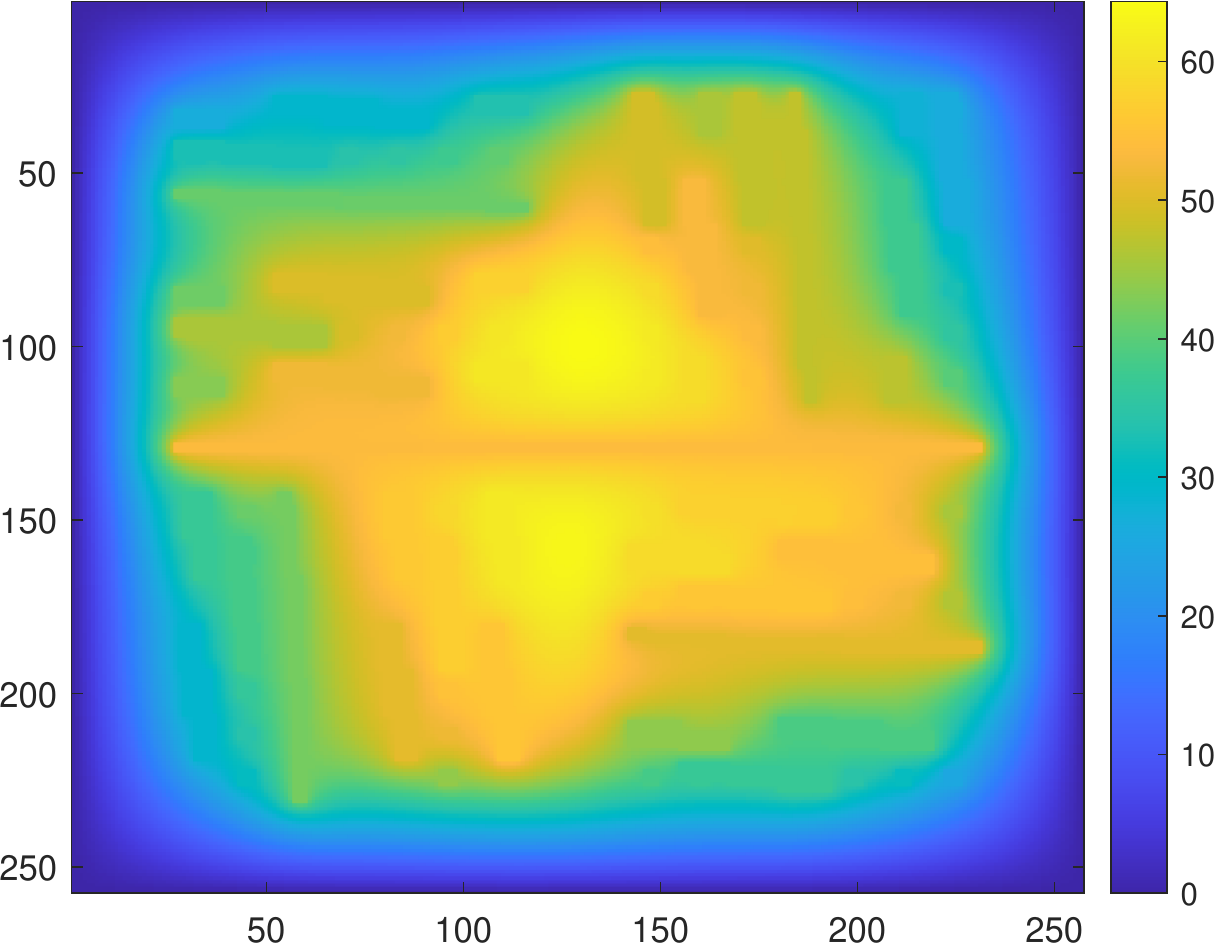}
	\end{subfigure}
	\begin{subfigure}{.45\textwidth}
		\centering
		\includegraphics[width = .9\linewidth]{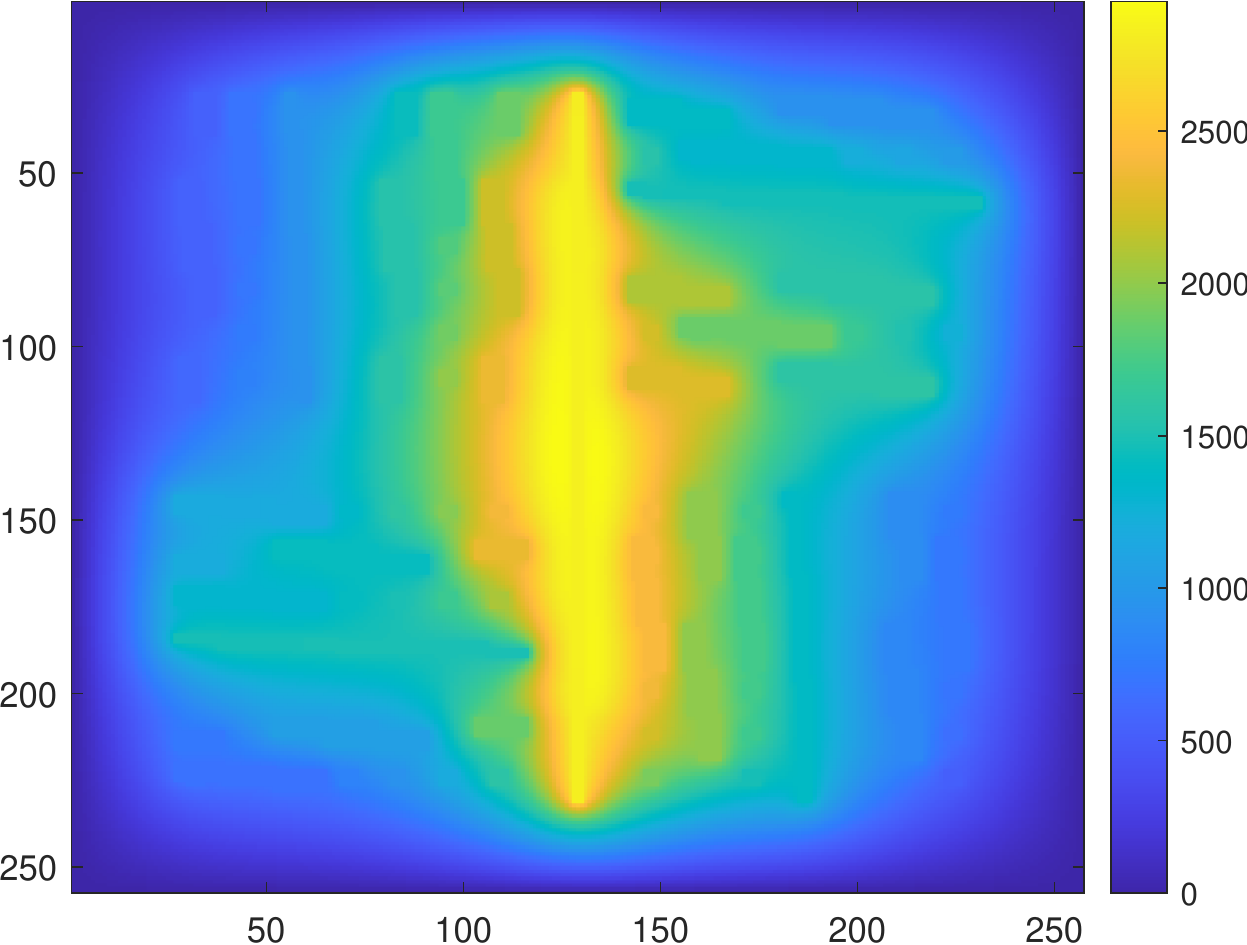}
	\end{subfigure}

	\begin{subfigure}{.45\textwidth}
		\centering
		\includegraphics[width = .9\linewidth]{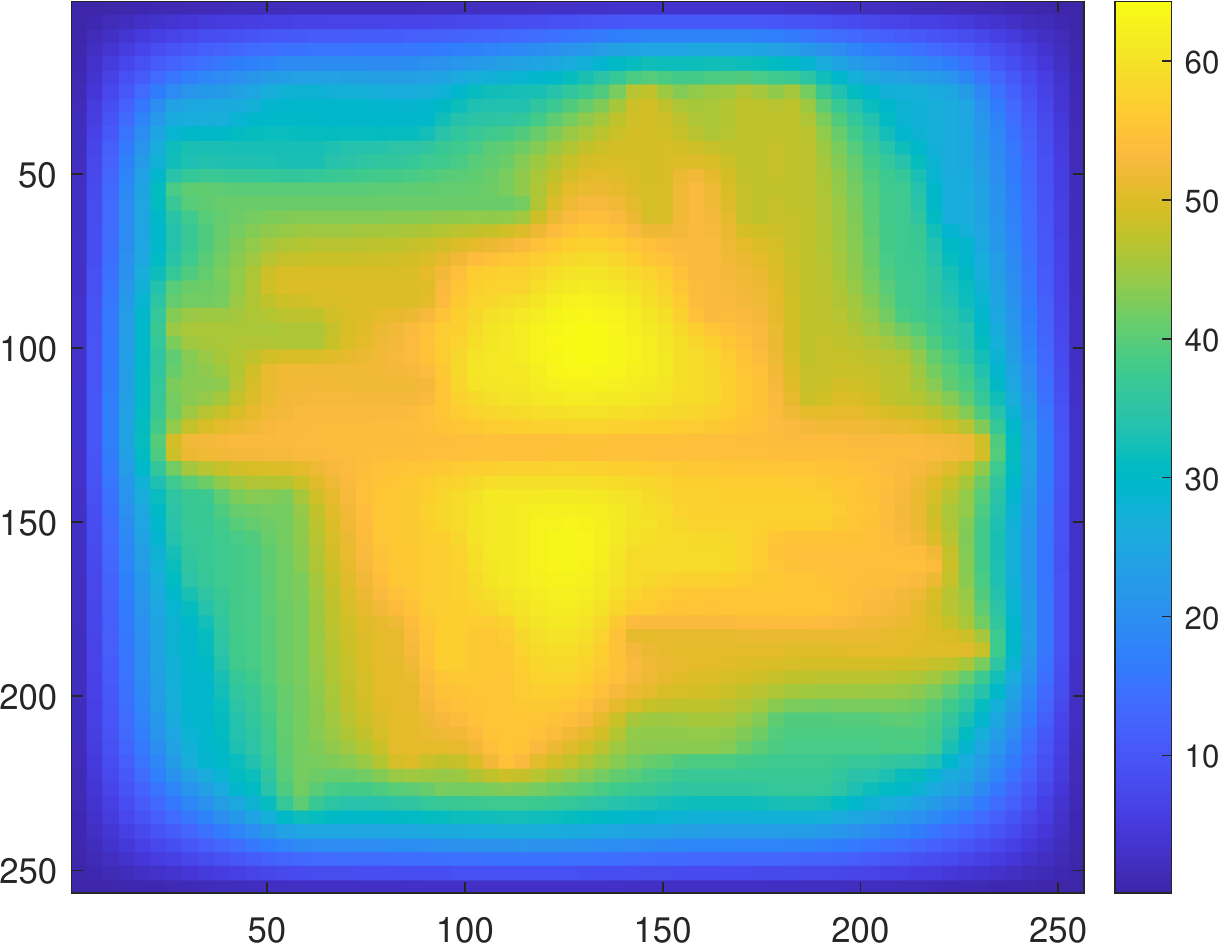}
	\end{subfigure}
	\begin{subfigure}{.45\textwidth}
		\centering
		\includegraphics[width = .9\linewidth]{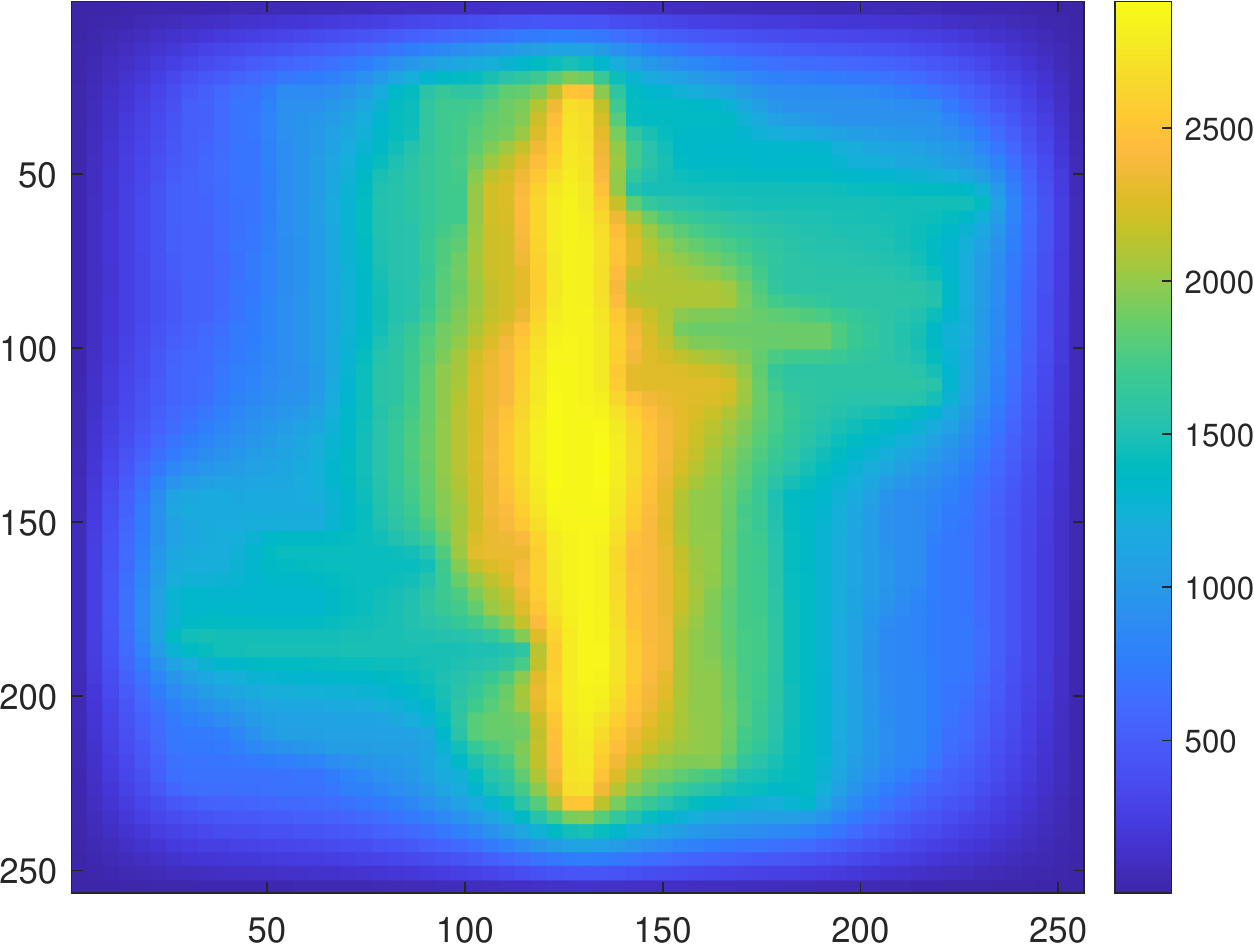}
	\end{subfigure}
	
	\begin{subfigure}{.45\textwidth}
		\centering
		\includegraphics[width = .9\linewidth]{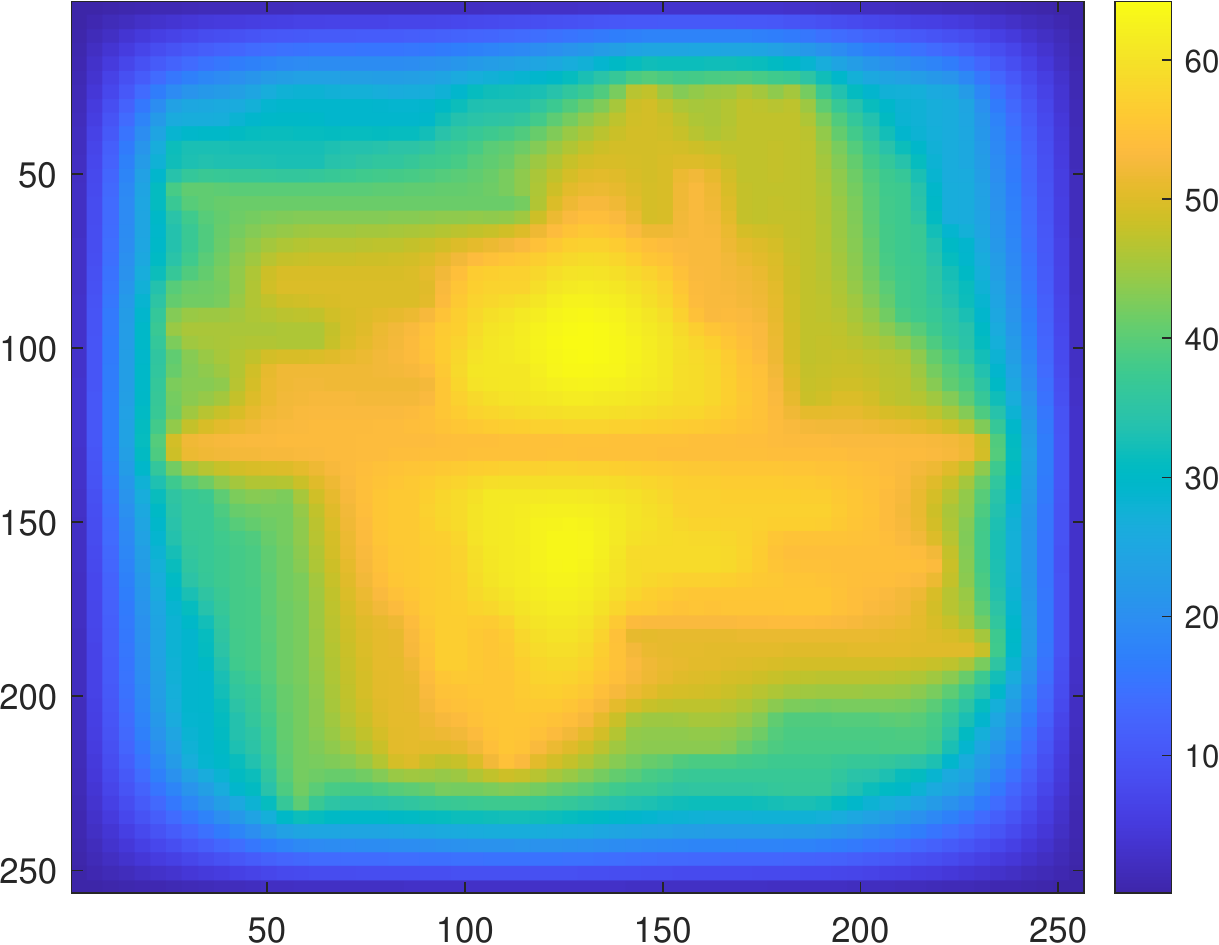}
	\end{subfigure}
	\begin{subfigure}{.45\textwidth}
		\centering
		\includegraphics[width = .9\linewidth]{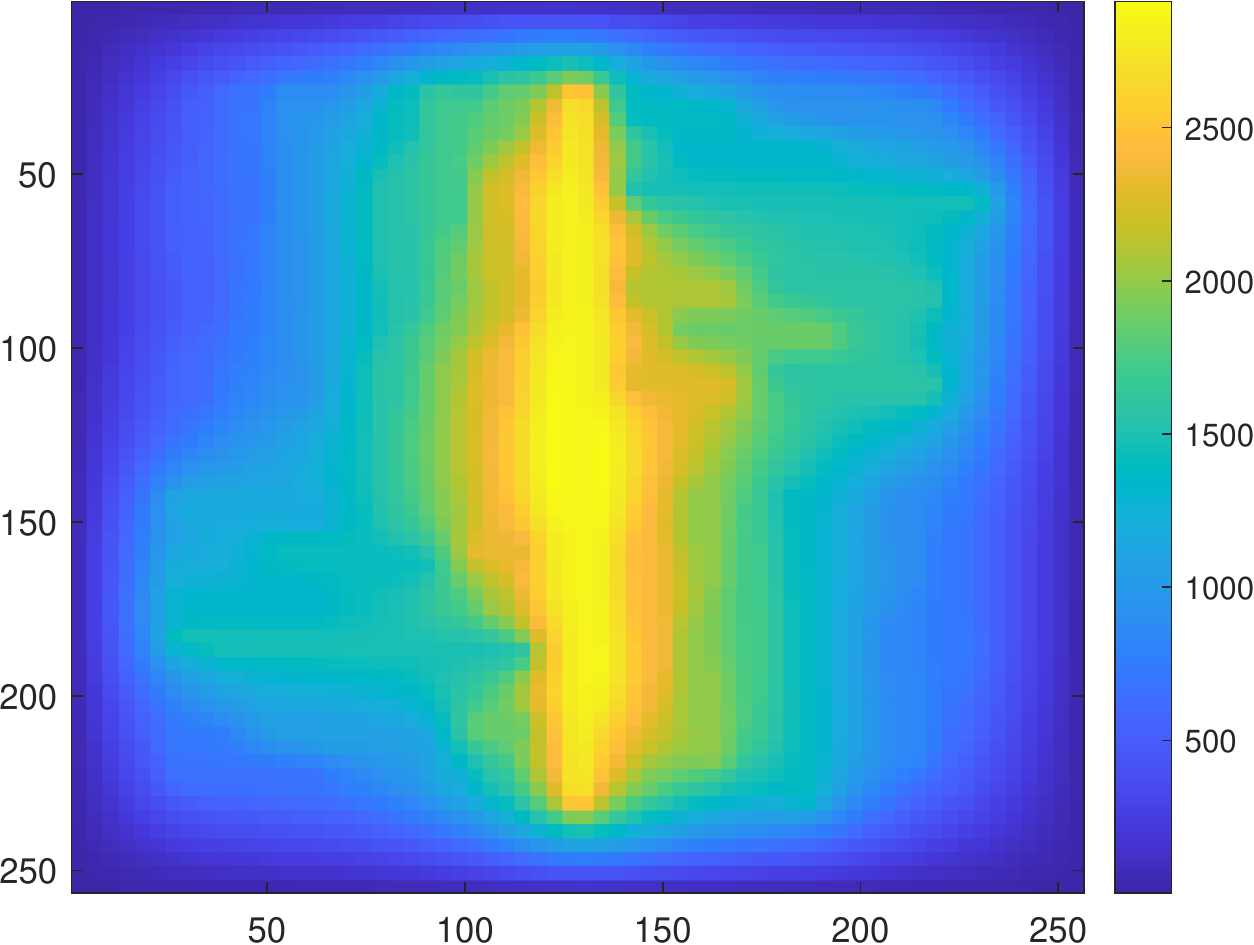}
	\end{subfigure}

	\caption{Plots of the numerical approximations of pressure with coarse mesh size $H=1/64$ and $m = 8$ oversampling layers in Experiment 1. 
	Left: first continuum. Right: second continuum. First row: fine-scale solution. Second row: coarse-scale average of fine-scale solution. Third row: NLMC solution.}
	\label{fig:sol_static}
\end{figure}

\begin{table}[!ht]
	\centering
	\begin{tabular}{c |c| c| c}
		\hline
		$H$ & $m$ & $e_{L^2}^{(1)}$ & $e_{L^2}^{(2)}$\\
		\hline
		1/8 & 3 & 96.7318\% & 88.9818\% \\
		1/16 & 5 & 32.3229\% & 19.3473\% \\
		1/32 & 6 & 0.6045\% & 0.3680\%\\
		1/64 & 8 & 0.0550\% & 0.0296\% \\
		\hline
	\end{tabular}
	\caption{Convergence of $e_{L_2}$ with respect to coarse mesh size $H$ in Experiment 1.}
	\label{tab:conv_static}
\end{table}

\begin{table}[!ht]
	\centering
	\begin{tabular}{c|c|c|c}
		\hline
		$m$ & Area Ratio & $e_{L^2}^{(1)}$ & $e_{L^2}^{(2)}$\\
		\hline
		3 & 4.79\% & 99.3677\% & 93.9357\% \\
		4 & 7.91\% & 76.6083\% & 55.1631\% \\
		5 & 11.81\% & 8.6605\% & 5.3115\% \\
		6 & 16.50\% & 0.6045\% & 0.3680\% \\
		\hline
	\end{tabular}
	\caption{Comparison of $e_{L_2}$ error with different number of oversampling layers $m$ for $H=1/32$ in Experiment 1.}
	\label{tab:error_32}
\end{table}
\begin{table}[!ht]
	\centering
	\begin{tabular}{c|c|c|c}
		\hline
		$m$ & Area Ratio & $e_{L^2}^{(1)}$ & $e_{L^2}^{(2)}$\\
		\hline
		2 & 0.61\% & 99.9102\% & 97.9631\% \\
		4 & 1.98\% & 99.1268\% & 91.9240\% \\
		6 & 4.13\% & 11.8898\% & 6.3181\% \\
		7 & 5.49\% & 0.7959\% & 0.4219\% \\
		8 & 7.06\% & 0.0550\% & 0.0296\% \\
		\hline
	\end{tabular}
	\caption{Comparison of $e_{L_2}$ error with different number of oversampling layers $m$ for $H=1/64$ in Experiment 1.} 
	\label{tab:error_64}
\end{table}

\subsection{Experiment 2: time-dependent case}
The time dependent case faces the similar issue with the error. $f_1(x,y) = 1$ and $f_2(x,y)$ is 
depicted in Figure~\ref{fig:source_time}, which represents a simplified five-spot well rate. 
The temporal domain is $[0, T]$ with final time $T= 5$.
In Figure~\ref{fig:sol_time}, we plot the fine-scale solution, the coarse-scale average 
and the NLMC coarse-scale solution with coarse mesh size $H = 1/64$ and number of oversampling layers $m = 8$. 
Again, the NLMC solution is a good approximation for the coarse-scale average. 
In Figure~\ref{fig:evolve_time}, we depict the change of pressure at different time steps. 
In Table~\ref{tab:conv_time}, we present the relative $L^2$ error with varying coarse grid size. 
Again, with the number of oversampling layers satisfying the sufficient condition, 
we can see that the error converges very well. 

\begin{figure}[!ht]
	\centering
	\includegraphics[width = .45\linewidth]{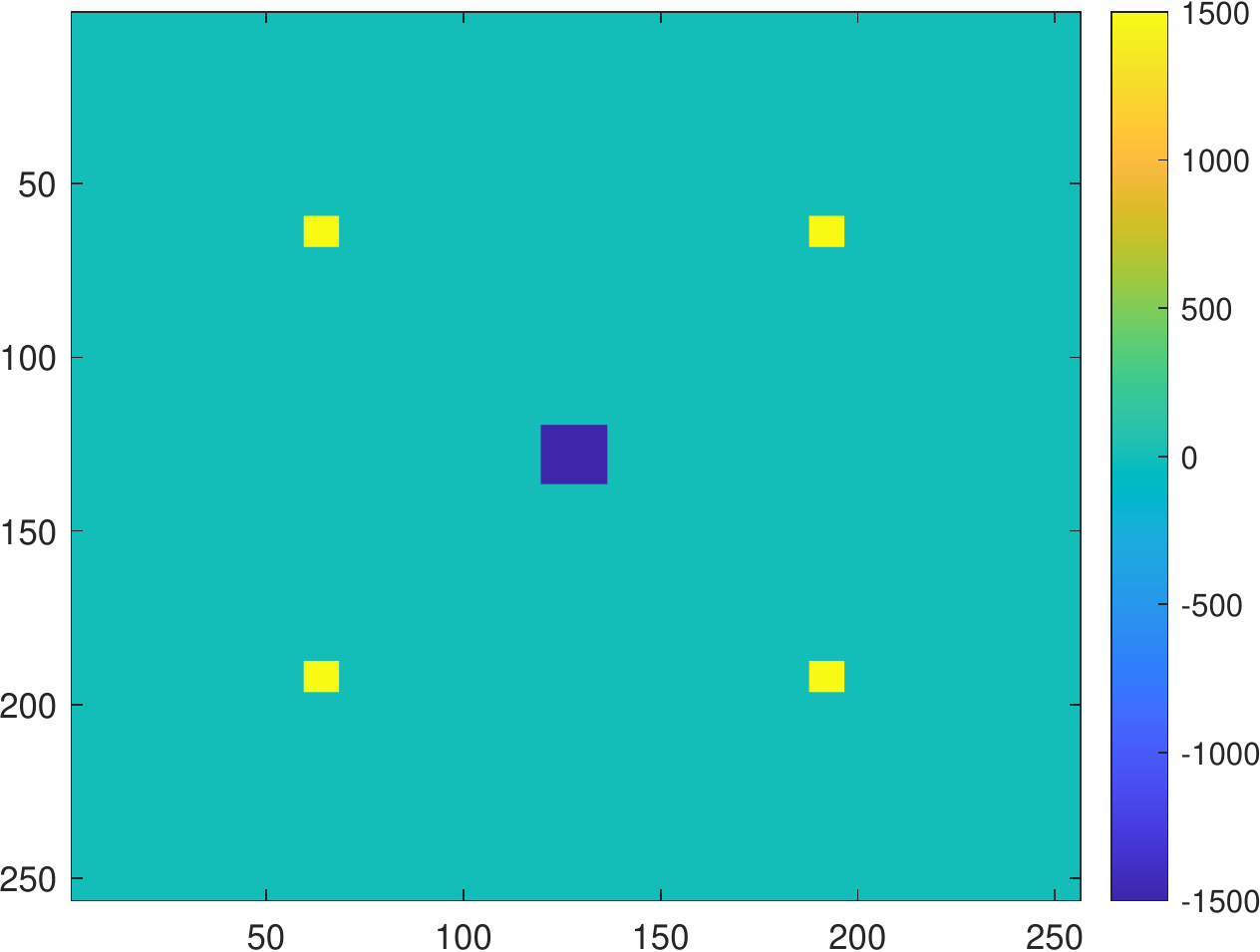}
	\caption{Source term $f_2$ in Experiment 2.}
	\label{fig:source_time}
\end{figure}

\begin{figure}[!ht]
	\centering
	\begin{subfigure}{.45\textwidth}
		\centering
		\includegraphics[width = .9\linewidth]{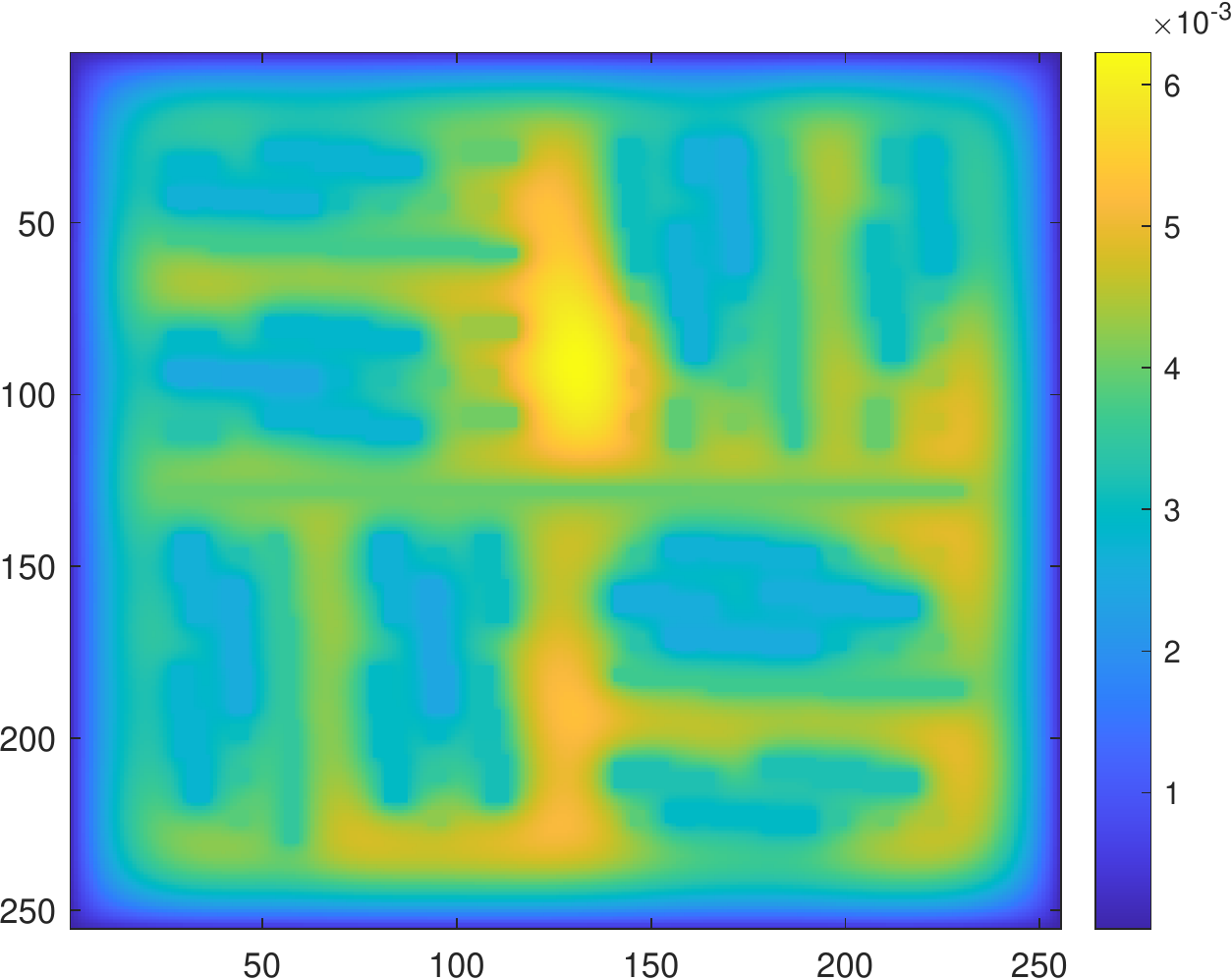}
	\end{subfigure}
	\begin{subfigure}{.45\textwidth}
		\centering
		\includegraphics[width = .9\linewidth]{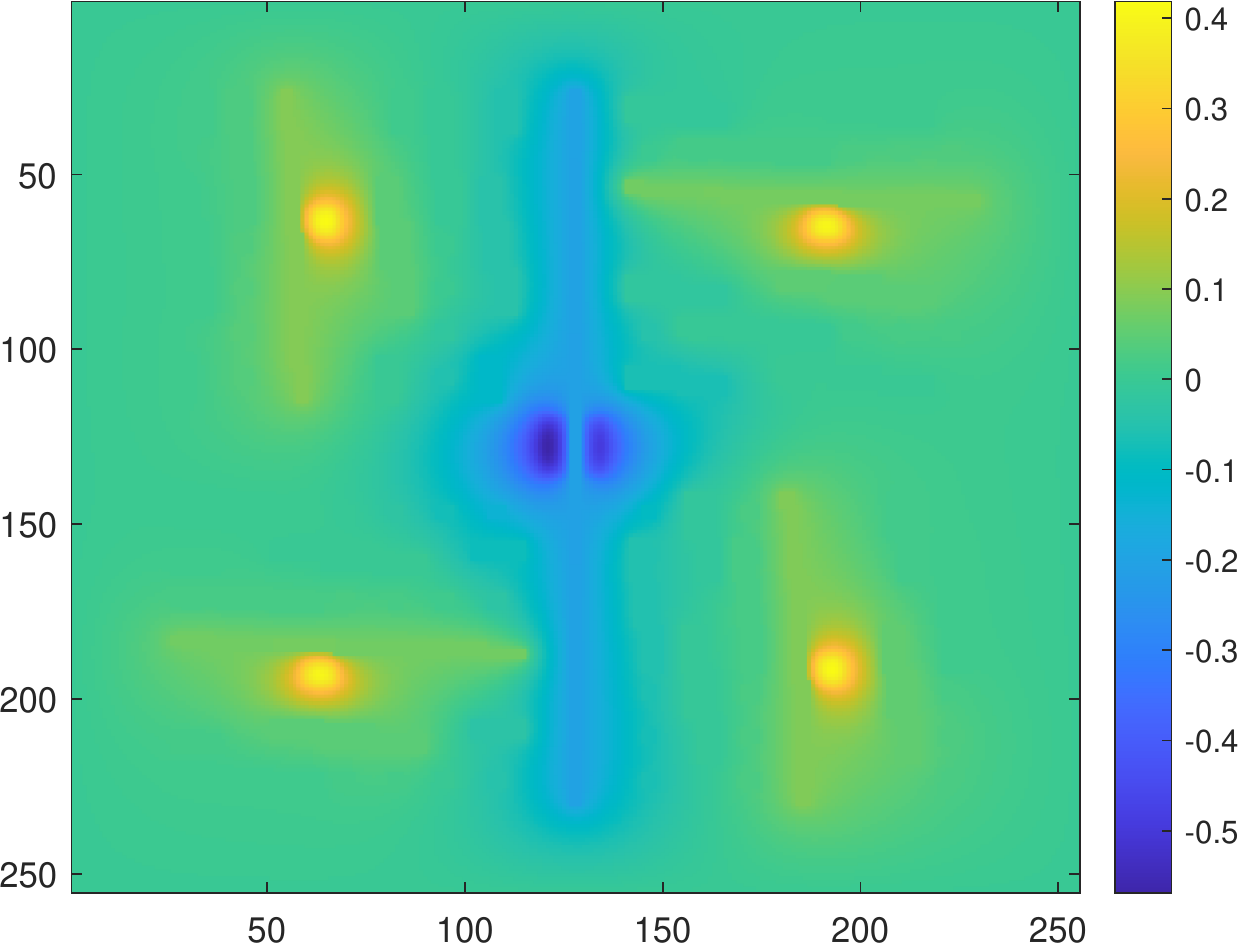}
	\end{subfigure}
	
	\begin{subfigure}{.45\textwidth}
		\centering
		\includegraphics[width = .9\linewidth]{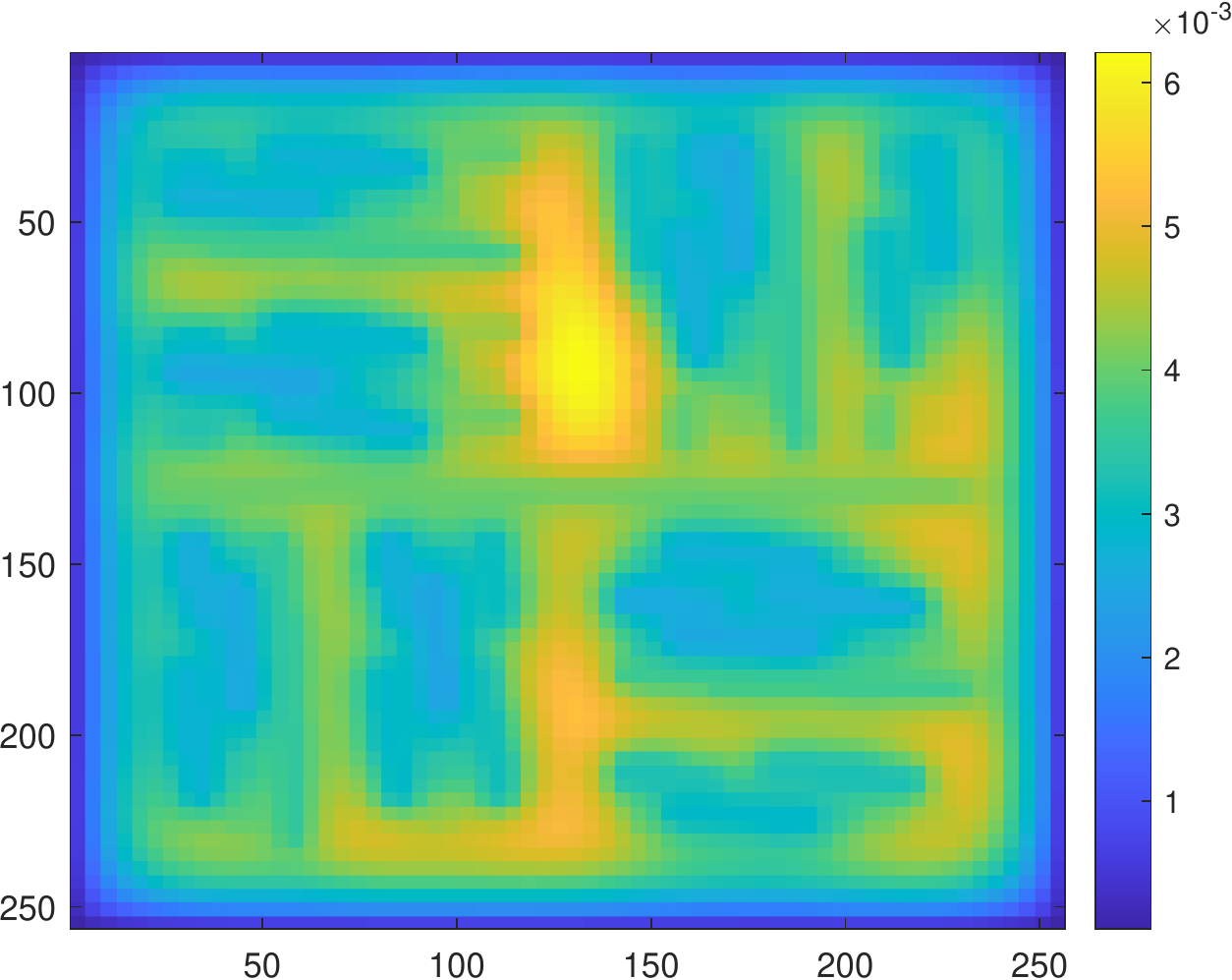}
	\end{subfigure}
	\begin{subfigure}{.45\textwidth}
		\centering
		\includegraphics[width = .9\linewidth]{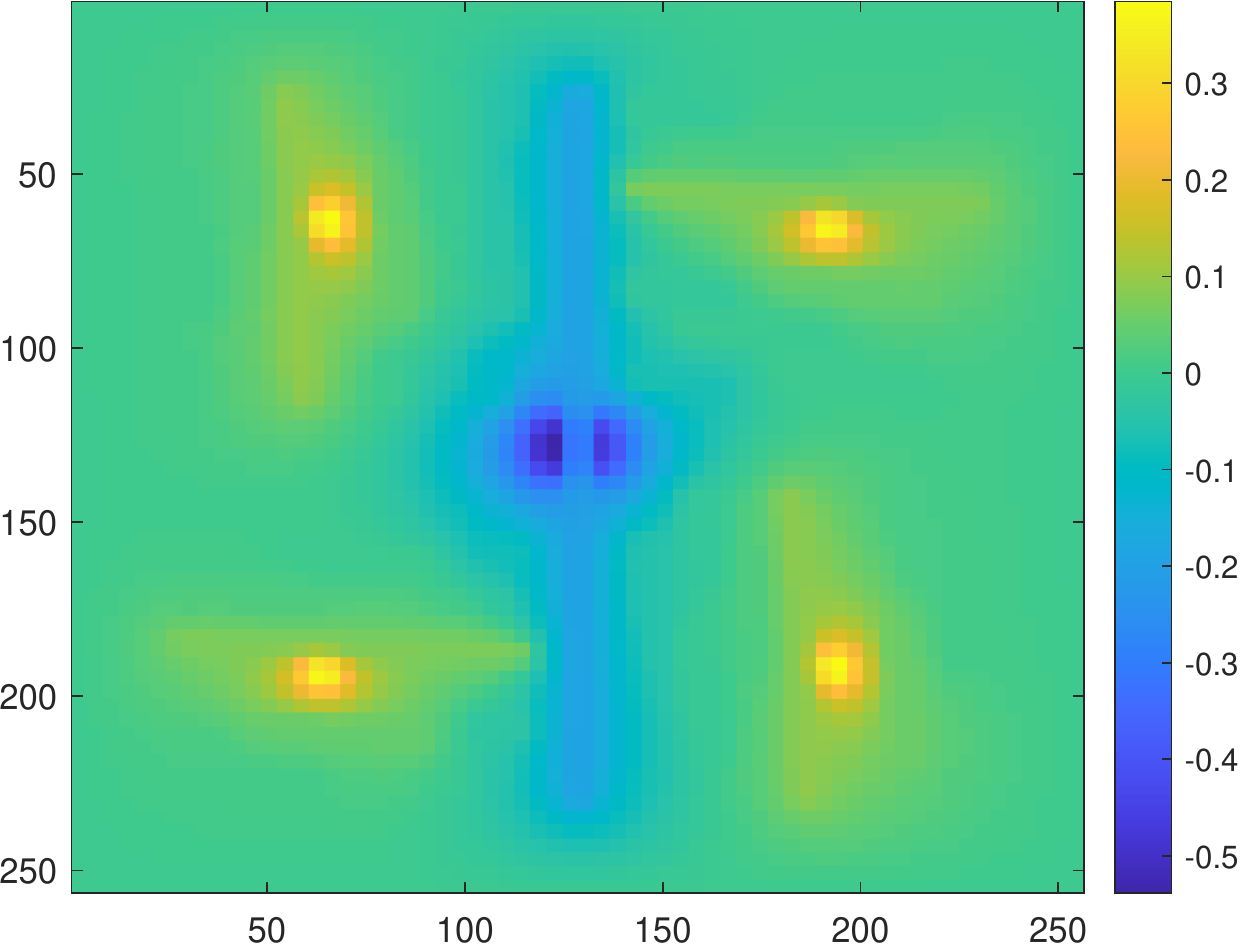}
	\end{subfigure}
	
	\begin{subfigure}{.45\textwidth}
		\centering
		\includegraphics[width = .9\linewidth]{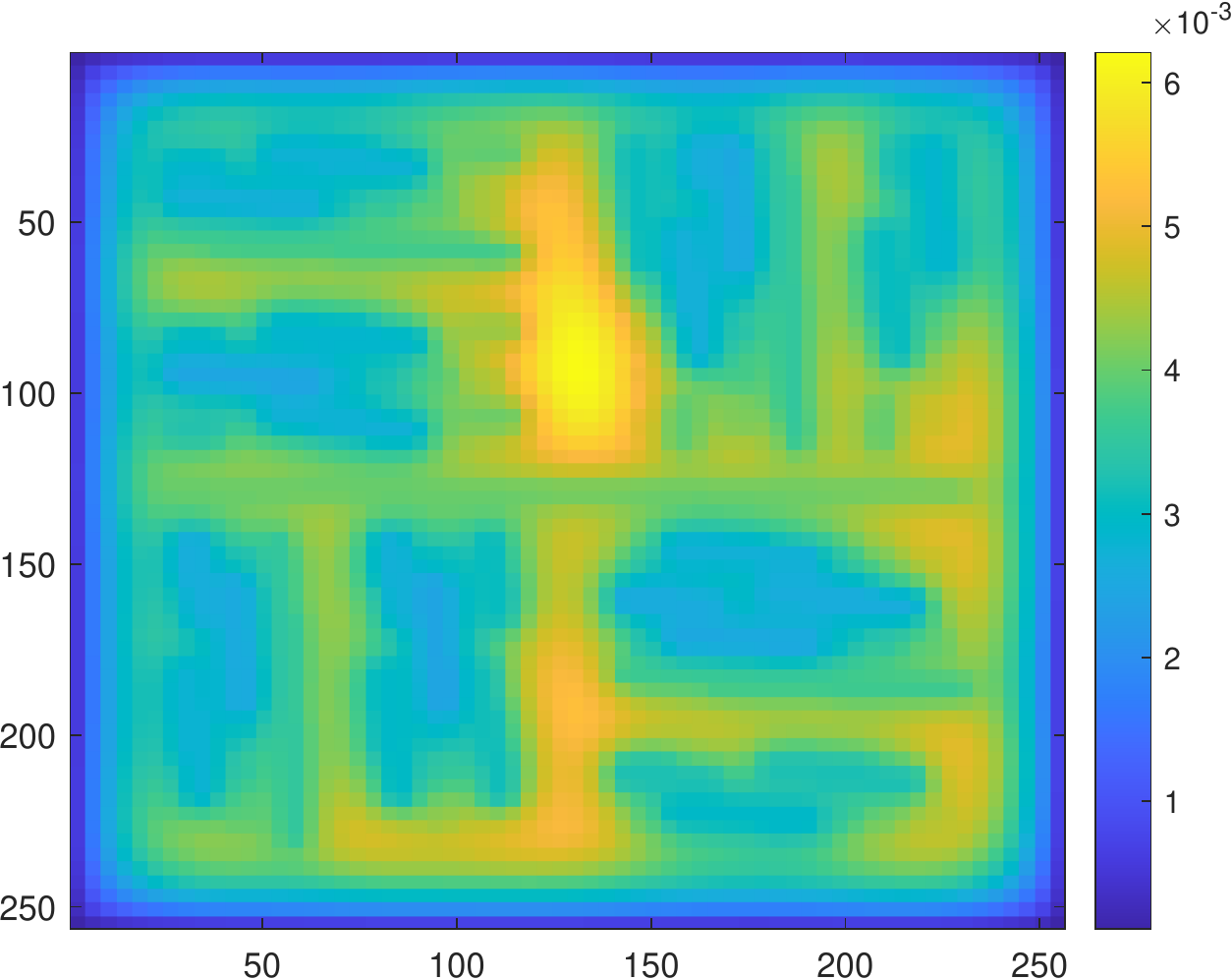}
	\end{subfigure}
	\begin{subfigure}{.45\textwidth}
		\centering
		\includegraphics[width = .9\linewidth]{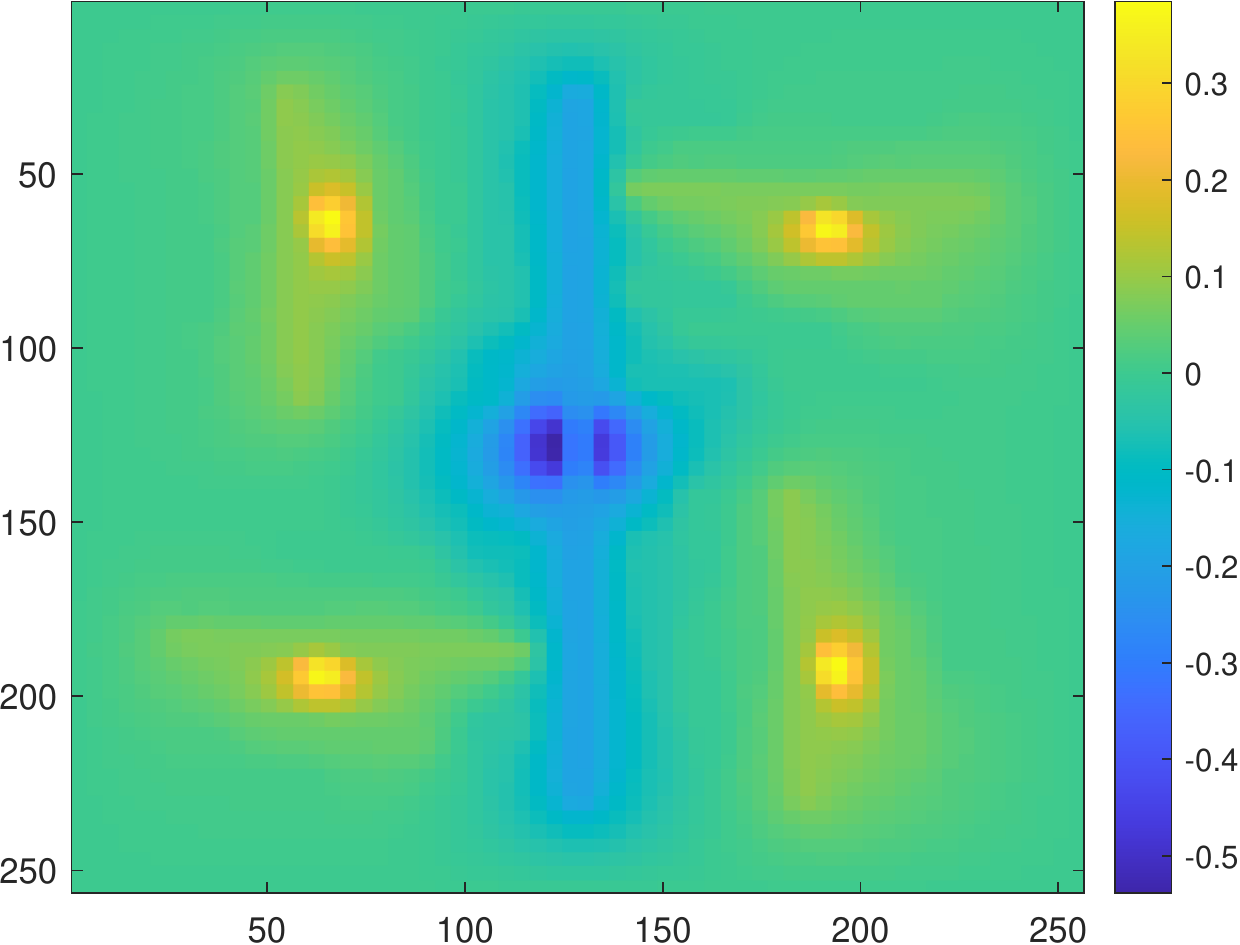}
	\end{subfigure}
	\caption{Plots of the numerical approximations of final-time pressure with coarse mesh size $H=1/64$ and $m = 8$ oversampling layers in Experiment 2. 
	Left: first continuum. Right: second continuum. First row: fine-scale solution. Second row: coarse-scale average of fine-scale solution. Third row: NLMC solution.}
	\label{fig:sol_time}
\end{figure}

\begin{table}[!ht]
	\centering
	\begin{tabular}{c |c| c| c}
		\hline
		$H$ & $m$ & $e_{L^2}^{(1)}$ & $e_{L^2}^{(2)}$ \\
		\hline
		1/8 & 3 & 20.6422\% & 58.4975\%\\
		1/16 & 5 & 1.1245\% & 2.6226\%\\
		1/32 & 6 & 0.0254\% & 0.0717\%\\
		1/64 & 8 & 0.0017\% & 0.0037\%\\
		\hline
	\end{tabular}
	\caption{Convergence of $e_{L_2}$ with respect to coarse mesh size $H$ in Experiment 2.}
	\label{tab:conv_time}
\end{table}

\begin{figure}[!ht]
	\centering
	\begin{subfigure}{.45\textwidth}
		\centering
		\includegraphics[width = .9\linewidth]{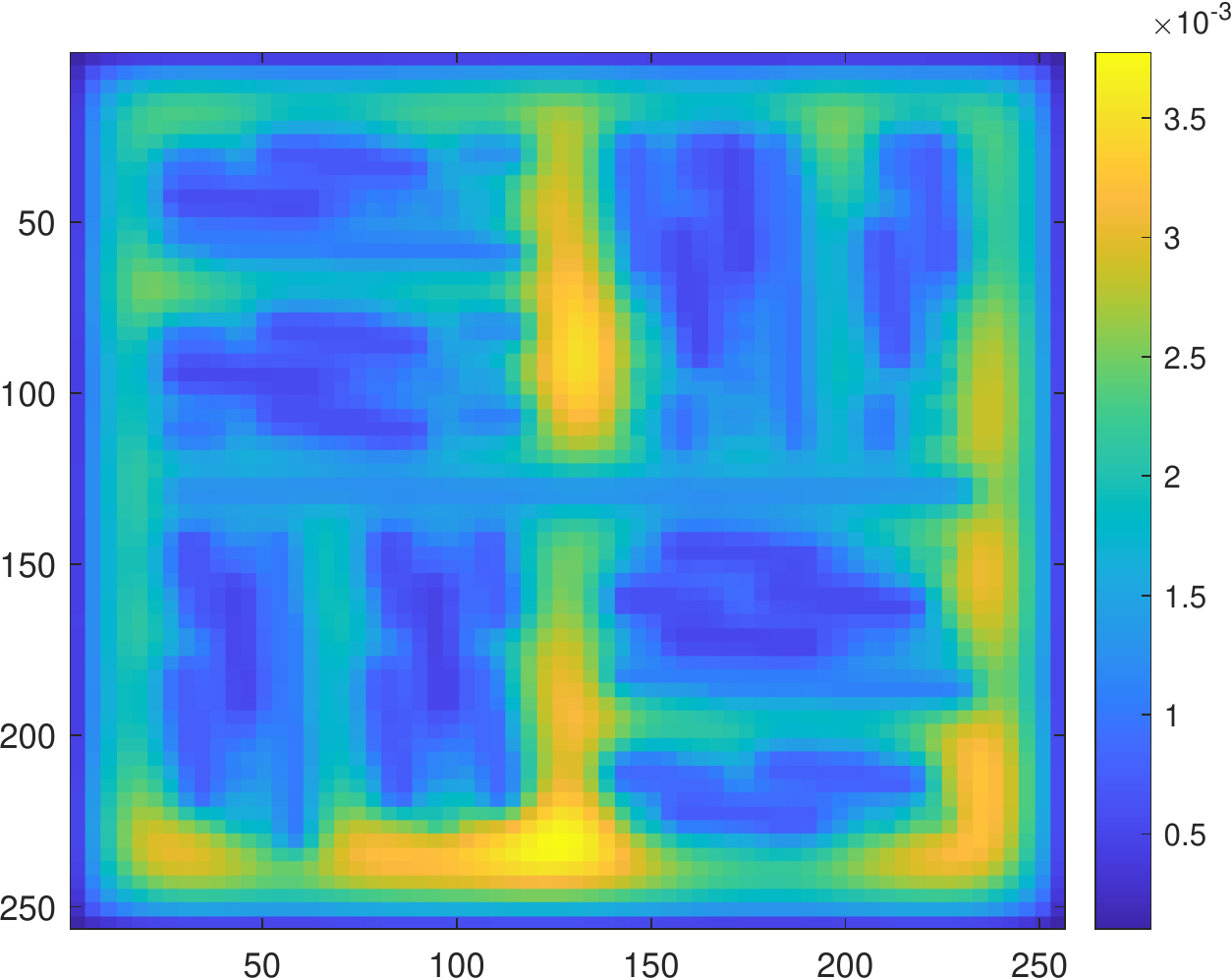}
	\end{subfigure}
	\begin{subfigure}{.45\textwidth}
		\centering
		\includegraphics[width = .9\linewidth]{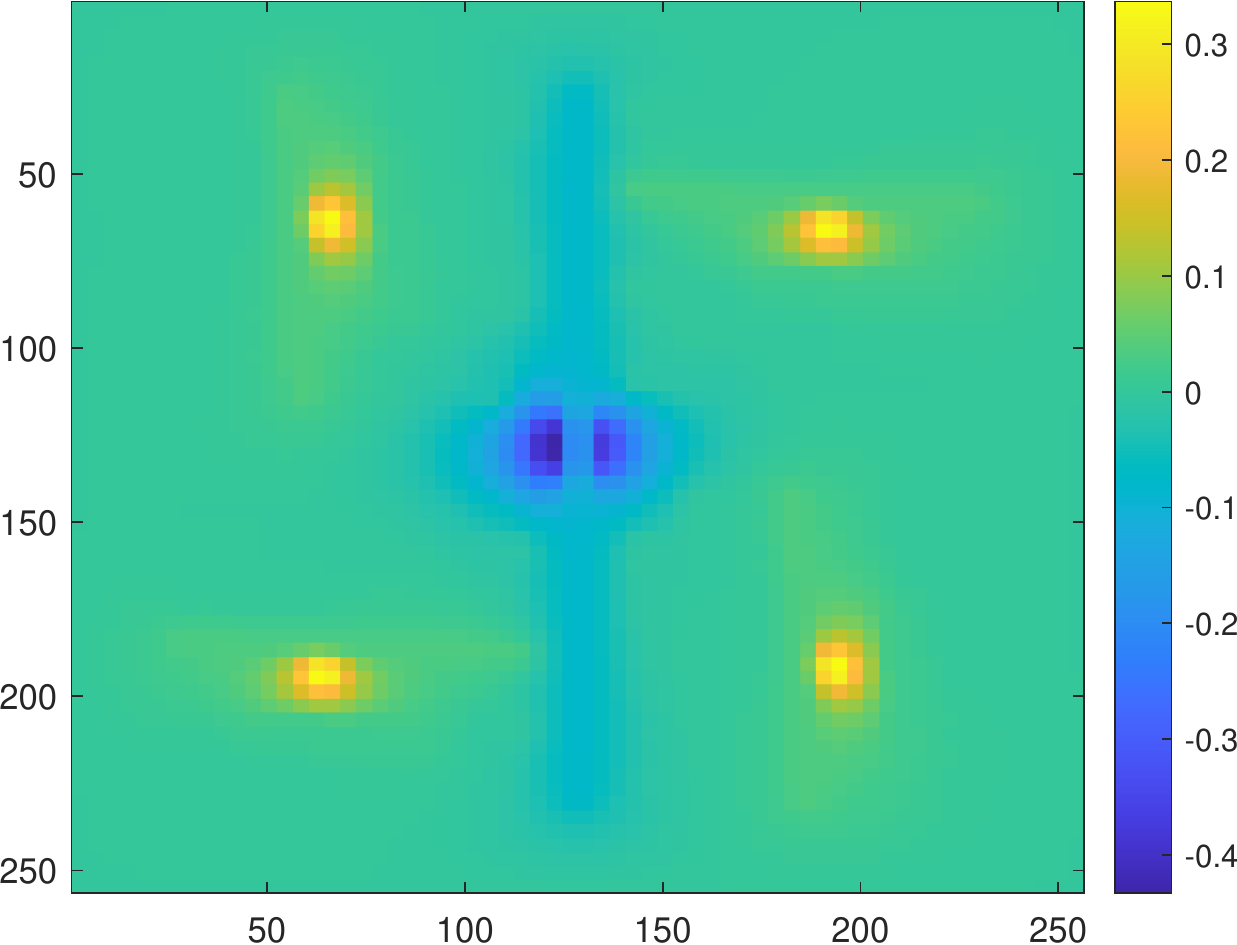}
	\end{subfigure}
	
	\begin{subfigure}{.45\textwidth}
		\centering
		\includegraphics[width = .9\linewidth]{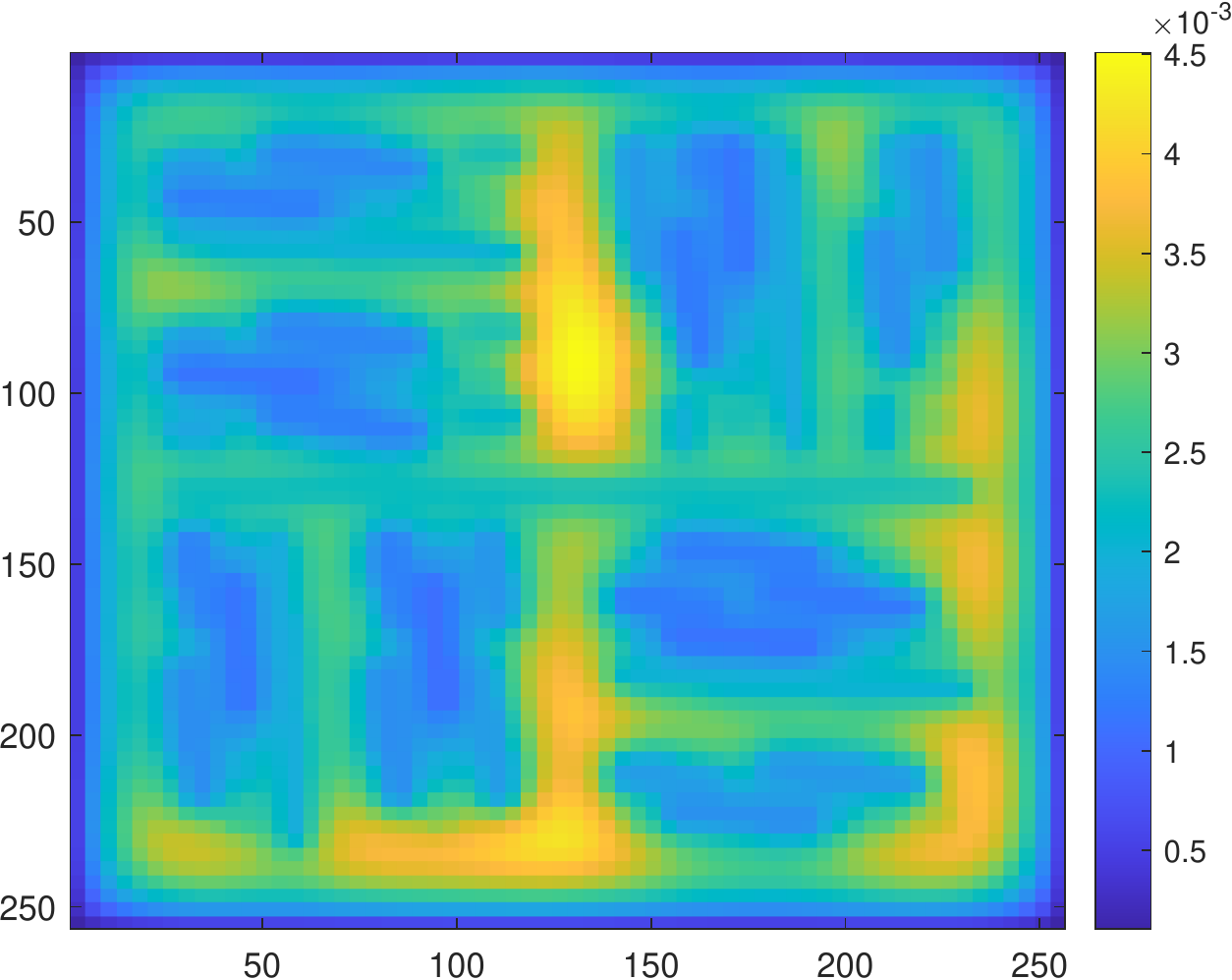}
	\end{subfigure}
	\begin{subfigure}{.45\textwidth}
		\centering
		\includegraphics[width = .9\linewidth]{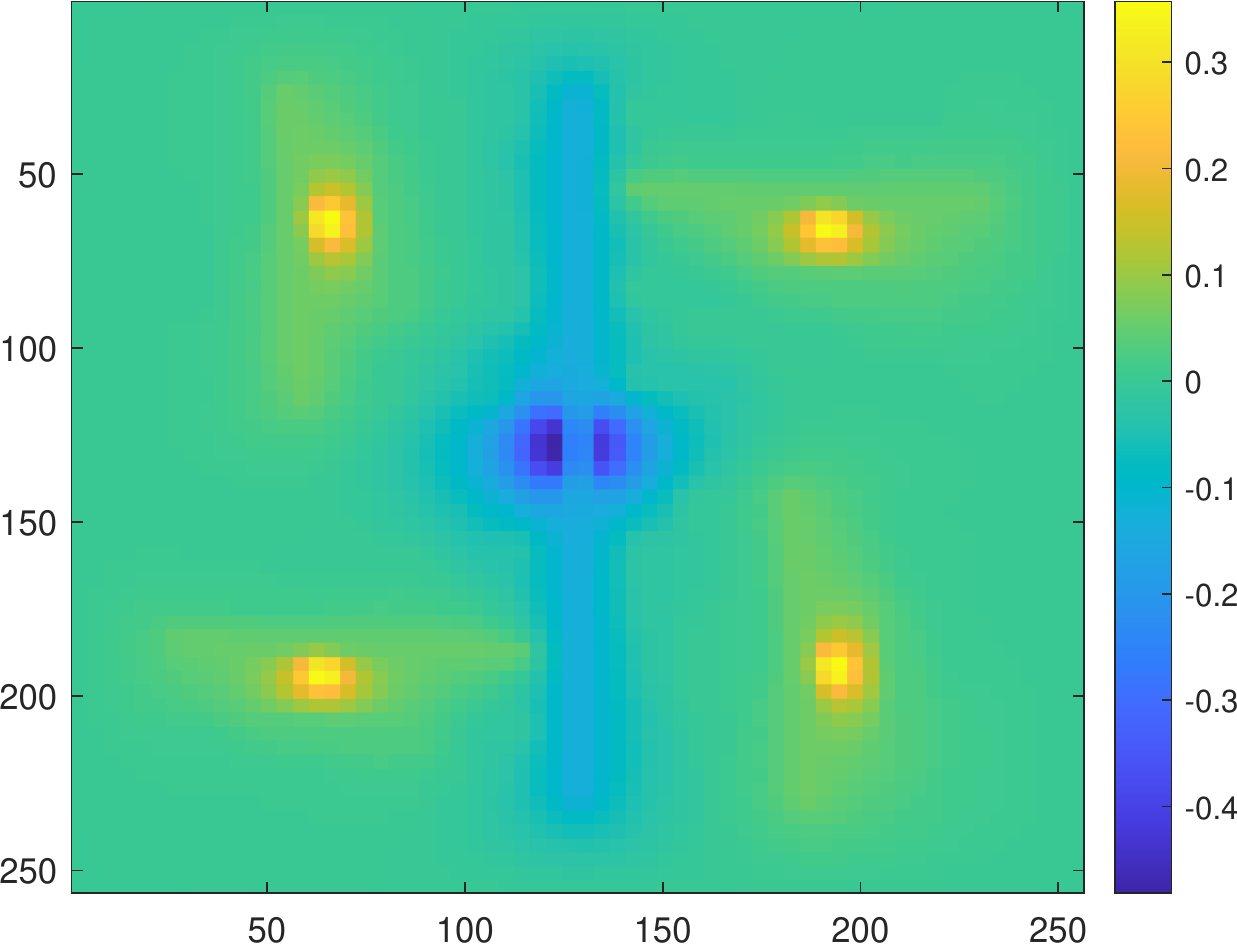}
	\end{subfigure}
	
	\begin{subfigure}{.45\textwidth}
		\centering
		\includegraphics[width = .9\linewidth]{time_nlmc_coarse_1_20}
	\end{subfigure}
	\begin{subfigure}{.45\textwidth}
		\centering
		\includegraphics[width = .9\linewidth]{time_nlmc_coarse_2_20}
	\end{subfigure}
	\caption{Plots of the NLMC numerical approximations of pressure at various time instants with coarse mesh size 
	$H=1/64$ and $m = 8$ oversampling layers in Experiment 2. 
	Left: first continuum. Right: second continuum. First row: $t = 1.25$. Second row: $t = 2.5$. Third row: $t = 5$.}
	\label{fig:evolve_time}
\end{figure}



\section{Conclusions}\label{sec:conclusions}

In this paper, we investigated the non-local multicontinuum (NLMC) upscaling method for a dual continuum model in fractured porous media. 
Localized multiscale basis functions that separate each continuum are constructed. 
To find the basis, we solve local problems subject to energy minimization constraints in oversampling coarse regions. 
It is showed that the basis functions equip the method with coarse mesh convergence. 
Some numerical examples are presented to support the theory. 
The numerical examples also indicate that the proposed method provides accurate and efficient coarse-grid approximation.

%

\vspace{6pt}


\bibliographystyle{plain}

\end{document}